\documentclass[11pt,a4paper]{article}

\usepackage{cite}
\usepackage{indentfirst,latexsym,bm}
\usepackage{mathrsfs}
\usepackage{CJK,fancybox,fancyhdr,color,graphicx,amsmath,amsthm,amssymb,enumerate,amsfonts}
\usepackage{caption}
\usepackage{algorithmicx,algorithm,algpseudocode}
\usepackage{anysize}
\usepackage{authblk}
\usepackage[misc]{ifsym}
\usepackage[T1]{fontenc}
\usepackage{subfigure}
\usepackage{xcolor}
\usepackage{colortbl,booktabs}

\marginsize{2cm}{2cm}{2.54cm}{2cm} 

\newtheorem{thm}{Theorem}[section]
\newtheorem{cor}[thm]{Corollary}
\newtheorem{lem}[thm]{Lemma}

\numberwithin{equation}{section}            
\renewcommand{\proof}{\noindent\textbf{Proof.}}

\title{Error estimates of Kaczmarz and randomized Kaczmarz methods}%

\author[a]{Chuan-gang Kang\thanks{Corresponding author. \authorcr \quad {\it E-mail address}: ckangtj@tjpu.edu.cn(C.-g. Kang)}} 
\author[a]{Heng Zhou}
\affil[a]{School of Mathematical Sciences, Tian Jin polytechnic University,Tian Jin 300387, China}
\date{}

\begin{document}

\maketitle


\begin{center}
\begin{minipage}{145mm}
{\bf Abstract:}  The Kaczmarz method is an iterative projection scheme for solving consistent system $Ax=b$. It is later extended to the inconsistent and ill-posed linear problems. But the classical Kaczmarz method is sensitive to the correlation of the adjacent equations. In order to reduce the impact of correlation on the convergence rate, the randomized Kaczmarz method and randomized block Kaczmarz method are proposed, respectively. In the current literature, the error estimate results of these methods are established based on the error $\|x_k-x_*\|_2$, where $x_*$ is the solution of linear system $Ax=b$. In this paper, we extend the present error estimates of the Kaczmarz and randomized Kaczmarz methods on the basis of the convergence theorem of Kunio Tanabe, and obtain some general results about the error $\|x_k-P_{N(A)}x_0-x^\dagger\|_2$.
\\ [5pt]
\mbox{\bf Key Words:}  Kaczmarz method, randomized Kaczmarz method, error estimates, convergent rate.
\\ [5pt]
\end{minipage}
\end{center}
\section{Introduction}

Kaczmarz method is a popular iterative method for solving linear systems. It was originally discovered by Kaczmarz in 1937. In 1970, it was rediscovered by Gordon, Bender and Herman\cite{Gordon1970} to deal with the computed tomography (CT). In the early stage of CT development, the computed techniques are dominated by back projection(BP)\cite{Long2010} and filtered back projection(FBP)\cite{Devaney1982A,Pan2004A} techniques because of the hardware constraints. With the development of hardware techniques, Kaczmarz method was widely used in CT reconstruction\cite{Herman1978,Natterer1986} for its outstanding properties, such as anti-interference and image reconstruction ability from the incomplete data\cite{Heffernan1983,Inouye1979,Reeds1987,Candes2006}

In 1971, Kunio Tanabe established the convergence theory\cite{Tanabe_1971} of Kaczmarz method for solving consistent linear systems and employed this method to solve the generalized inverse of singular matrices.  Moreover, he verified that Kaczmarz method works well for both singular and non-singular systems.

But the issue of the convergent rate is very challenging. The convergent rate of the method depends strongly on the ordering of the equations. Subsequently, using the rows of matrix $A$ in Kaczmarz method in a random order rather than original order, which can often substantially improve the convergence, was proposed as a computed techniques and was named randomized Kaczmarz method\cite{Herman_1993} later. The random selection of the rows from coefficient matrix $A$ weaken the influence of the original order of the equations and are in line with the actual situation, so randomized Kaczmarz method is appealing for practical applications.

In 2009, Strohmer and Vershynin established the exponential convergence\cite{Strohmer_2009} of randomized Kaczmarz method (RKM), where the convergent rate depends on a variant of the condition number, for consistent and inconsistent linear systems.

In 2014 Needell and Tropp considered a block algorithm\cite{Needell2014paved} that used a randomized control scheme to choose the subspace at each step.  They analyzed the convergent rate of randomized block Kaczmarz method\cite{Needell2014paved} for overdetermined least-square problems in that paper and obtained some important and meaningful results.
In 2015, Needell\cite{Needell2015Randomized} and his collaborators analyzed two block versions of Kaczmarz method each with a randomised projection, designed to converge in expecatation to the least squares solution.

For noisy linear systems, Elfing, Hansen and Nikazad proved the semi-convergent behavior of Kaczmarz method\cite{Elfing_2014} for solving singular linear systems and ill-posed problems.  They illustrated that the Kaczmarz method is semi-convergent for ill-posed problem. The semi-convergent behavior actually clarify that Kaczmarz method can be considered as a kind of regularization method\cite{Hanke97Regularizing-properties,Engl89,Engl96}for solving inverse problems or ill-posed problems.

These convergent rate results do not fully explain the excellent empirical performance of randomized Kaczmarz method for solving linear inverse problems, especially in the case of noisy data. In fact, these convergent rate more reflect the change of the error with noise-free term rather than the perturbance, i.e., the decrease of the iterative error has nothing to do with the perturbed term.

Yuling Jiao, Bangti Jin and Xiliang Lu considered the relationship between the iterative error and perturbed term by splitting the iterative error into high-frequency error and low-frequency error\cite{Yuling_Jiao_2017}. They discovered the change process of high- and low-frequency error.

Most of the present literature about convergent rate depend on the hypothesis that $A$ has full column rank or $Ax=b$ is consistent.  Few literature considered the least square problems but they assume that $Ax=b$ has a unique solution, so they could not consider fully the property of generalized solution. Furthermore, the current literature ignore the influence of the inial value $x_0$ or simply make $x_0=0$.

Although $x_0$ can be fixed as zero, the reasonable case is that $x_0$ is arbitrarily chosen.  In this paper, inspired by the opinion of some literatures \cite{Tanabe_1971,Yuling_Jiao_2017,Ma2015Convergence,Needell2015Randomized}, we consider the convergent rate of Kaczmarz method, randomized Kaczmarz method and randomized block Kaczmarz method, and establish the convergent rate results on the basis of $\|x_k-P_{N(A)}x_0-x^\dagger\|_2$.

The rest of the paper is organized as follows. In section 2, we analyze the convergent rate of randomized Kaczmarz and classical Kaczmarz method, and illustrate the convergent behavior of these methods. Then in section 3, on the basis of block division $S=S_1\bigcup S_2\bigcup\cdots\bigcup S_r$ about index set $S=\{1,2,\cdots,m\}$, we propose a randomized Kaczmarz method and  analyze its convergent rate. Last, In section 4, we illustrate the convergent behavior through several classical numerical experiments.

\section{Kaczmarz and randomized Kaczmarz methods}
In this paper, the solution of the linear equations
\begin{align}\label{linear.system}
   Ax=b
\end{align}
are considered, where $A=(a_{ij})$ is an $m\times n$ complex matrix, $x$ and $b$ are $n-$ and $m-$dimension complex column vectors, respectively. we shall denote by $a_i=(a_{i1},\cdots,a_{in})^T,i=1,2,\cdots,m$ and $b_i$ the $i$-th row of $A$ and the $i$-th component of the vector $b$, respectively. We shall suppose that $a_i\neq 0, i=1,2,\cdots,m$.

In perturbed case, the right-hand side of \eqref{linear.system} is $b^\delta$, where $\|b^\delta-b\|_2\le \delta$. The solution of \eqref{linear.system} may or may not exist, furthermore, the least square solution may not be unique.

The classical Kaczmarz method\cite{Kaczmarz1937} can be described as
\begin{subequations}
\begin{align}\label{kaczmarz.1}
  x_k=x_{k-1}+\frac{b_i-(a_i,x_{k-1})}{\|a_i\|_2^2}a_i,k=1,2,\cdots,
\end{align}
\text{or}
\begin{align}\label{kaczmarz.2}
  x_k=(I-\frac{a_ia_i^T}{\|a_i\|_2^2})x_{k-1}+\frac{b_i}{\|a_i\|_2^2},k=1,2,\cdots,
\end{align}
\end{subequations}
where $i=(k\bmod m)+1$.

The randomized Kaczmard method\cite{Herman_1993,Natterer1986} can be described as
\begin{align} \label{randomized_kaczmarz}
  x_k=x_{k-1}+\frac{b_i-(a_i,x_{k-1})}{\|a_i\|_2^2}a_i,k=1,2,\cdots,
\end{align}
where $i$ is chosen in the set $\{1,2,\cdots,n\}$ by the probability $p_r=\frac{\|a_i\|_2^2}{\|A\|_F^2}$.

The convergence of Kaczmarz method was obtained by Tanabe\cite{Tanabe_1971} as the following theorem.
\begin{thm} \label{thm.1}
     \normalfont For any $m\times n$ matrix $A$ with nonzero rows and any $m$-dimensional column vector $b$, the algorithm \eqref{kaczmarz.1}
         generates a convergent sequence of vectors such that
    \begin{align}\label{tanabe.res}
    \lim\limits_{x\rightarrow \infty}x_k=P_{N(A)}x_0+Gb.
  \end{align}
where $x_0\in \mathcal{R}^n$ is an arbitrary initial vector.
\end{thm}

In \eqref{kaczmarz.2}, if let
\begin{align*}
    P_i=I-\frac{a_ia_i^T}{\|a_i\|_2^2}, ~~~~~~i=1,2,\cdots,m
\end{align*}
and
\begin{align*}
    Q_i=P_1P_2\cdots P_i,~~~~~~i=1,2,\cdots,m,
\end{align*}
where $Q_0=I, Q=Q_m$, and denote
\begin{align*}
  R=(\frac{Q_0a_1}{\|a_1\|_2^2}, \frac{Q_1a_2}{\|a_2\|_2^2},\cdots,\frac{Q_{m-1}a_m}{\|a_m\|_2^2}).
\end{align*}
Thereby, the operator $G$ in \eqref{tanabe.res} can be described as
\begin{align*}
    G=(I-\widetilde{Q})^{-1}R=\sum\limits_{j=0}^\infty Q^jR.
\end{align*}

In 2015, Kang and Zhou\cite{kang2015} analyzed the generalized inverse property of operator $G$ and proved that $Gb$ is not only the least square solution but also the minimal norm solution, i.e. $Gb$ is the Moore-Penrose generalized solution of $Ax=b$.

The convergent rate analysis in the following of this paper is mostly on the basis of this theorem.

The following theorem summarizes typical convergence results of randomized Kaczmarz method for consistent and inconsistent linear systems\cite{Needell_2010,Strohmer_2009,Zouzias_2013,Yuling_Jiao_2017}.

\begin{thm}\label{thm.randomized.Kaczmarz}
\normalfont Let $x_k$ be the solution generated by RKM at iteration $k$. and $\kappa_A=\|A\|_F^2\|A^\dagger\|_2$ be a (generalized) condition number. Then the following statements hold.

(i) For exact data, there holds
\begin{align}
E\left[\|x_k-x^*\|_2^2\right]\le (1-\kappa_A^{-2})^k\|x_0-x^*\|_2^2.
\end{align}

(ii) For noisy data, there holds
\begin{align}
E\left[\|x_k-x^*\|_2^2\right]\le (1-\kappa_A^{-2})^k\|x_0-x^*\|_2^2+\frac{\delta^2}{\sigma_{\min}^2(A)}.
\end{align}
where $x^*$ is a solution of the consistent system $Ax=b$.
\end{thm}

In \cite{Needell_2010,Strohmer_2009,Zouzias_2013}, The error estimate about $\|x_{k+1}-x^*\|_2$ was considered for consistent system or some inconsistent system. For any initial vector $x_0$, Kunio Tanabe has proved that the limit of Kacmarz method is $P_{N(A)}x_0+Gb$ for consistent linear system, where $Gb$ has been proved to be a Moore-Penrose generalized solution \cite{Ben-Israel_2003,kang2015,Elfing_2014}. Inspired by them, we will consider the theory analysis based on the error term $\|x_{k+1}-P_{N(A)}x_0-x^\dagger\|_2$ rather than $\|x_{k+1}-x^*\|_2$.

For any $x_0\in R^n$, define manifold $D_r=P_{N(A)}x_0+N(A)^\bot$ and restrict the operator $A$ on the manifold $D_r$, that is
\begin{align}
  A|_{D_r}:P_{N(A)}x_0+N(A)^\bot\mapsto R(A).
\end{align}
\begin{lem}\label{lemma.1}\normalfont
  For any $x_0\in R^n$, assume $D_r=P_{N(A)}x_0+N(A)^\bot$, then for the vector sequence $x_k$ generated from \eqref{randomized_kaczmarz}, there hold $x_k\in D_r, k=1,2,\cdots$.
\end{lem}
\proof ~we prove this lemma with mathematical induction.

First, for any initial vector $x_0\in R^n$, obviously there has $x_0\in D_r$. Next, we will prove $x_1\in D_r$, from randomized Kaczmarz iteration \eqref{randomized_kaczmarz}, there is
\begin{align}
  x_1=x_0+\frac{b_i-(a_i,x_0)}{\|a_i\|_2^2}x_i.
\end{align}
For any $\tilde{x} \in N(A)$,
\begin{align}
  (x_1-P_{N(A)}x_0,\tilde{x})&=( (I-P_{N(A)})x_0+\frac{b_i-(a_i,x_0)}{\|a_i\|_2^2}a_i,\tilde{x} )\nonumber\\
                             &=\big{(}(I-P_{N(A)}\big{)}x_0,\tilde{x})+(\frac{b_i-(a_i,x_0)}{\|a_i\|_2^2}a_i,\tilde{x})\nonumber\\
                             &=\frac{b_i-(a_i,x_0)}{\|a_i\|_2^2}a_i^T\tilde{x}=0\nonumber,
\end{align}
which shows that $x_1=P_{N(A)}x_0+N(A)^T\in D_r$.

Second, suppose $x_k\in D_r$, we prove $x_{k+1}\in D_r$ in the following. From the hypothesis $\tilde{x}\in N(A)$, so
\begin{align}
  (x_{k+1}-P_{N(A)}x_0,\tilde{x})&=(x_k-P_{N(A)}x_0+\frac{b_i-(a_i,x_k)}{\|a_i\|_2^2}a_i,\tilde{x})\nonumber\\
                             &=(x_k-P_{N(A)}x_0,\tilde{x} )+(\frac{b_i-(a_i,x_k)}{\|a_i\|_2^2}a_i,\tilde{x}).
\end{align}
Because of $x_k\in D_r$ and $\tilde{x}\in N(A)$, there holds $\left(x_k-P_{N(A)}x_0,\tilde{x} \right)=0$, therefore there has
\begin{align*}
  (x_{k+1}-P_{N(A)}x_0,\tilde{x})=(\frac{b_i-(a_i,x_k)}{\|a_i\|_2^2}a_i,\tilde{x})=0,
\end{align*}
which shows $x_{k+1}\in D_r$.

To sum up, for any $x_0\in R^n$, the vector sequence $\{x_k\}_{k=1}^\infty$ generated from randomized Kaczmarz method belong to the manifold $D_r$. \hfill{$\square$}

The Kaczmraz method is an orthogonal projection method, so there holds the follow lemma.
\begin{lem}\label{orthogonal}\normalfont
  Let $u,v\in R^n$ are arbitrary vectors, then the vector $(I-\frac{uu^T}{\|u\|_2^2})v$ and $u$ are orthogonal.
\end{lem}
\proof ~~
\begin{align*}
    \big{(}(I-\frac{uu^T}{\|u\|_2^2})v, u\big{)}=(v,u)-(\frac{uu^T}{\|u\|_2^2}v, u)=(v, u)-(v, \frac{uu^T}{\|u\|_2^2}u)=(v, u)- (v, u)=0.
\end{align*}

\begin{thm}\label{thm.5}
\normalfont
For the linear system $Ax=b$, where $A\in R^{m\times n}, b\in R^m$, Assume $A^\dagger$ is the Moore-Penrose inverse of $A$, $x^\dagger$ is the Moore-Penrose generalized solution of $Ax=b$, $P_{N(A)}:R^n\rightarrow N(A)$ and $Q: R^m\rightarrow R(A)$ are orthogonal projection operators, $\{x_k\}_{k=1}^\infty$ is the vector sequence generalized from randomized Kaczmarz method \eqref{randomized_kaczmarz}, $\kappa_A=\|A\|_F\|A^\dagger\|_2$ is generalized condition number, then
there holds
\begin{align}\label{thm.result.1}
  E\big{[}\|x_{k+1}-P_{N(A)}x_0-x^\dagger\|_2^2\big{]}\le (1-\frac{1}{\kappa^2(A)})^k\|x_0-P_{N(A)}x_0-x^\dagger\|_2^2+\frac{1}{\|A\|_F^2}\|(I-Q)b\|_2^2.
\end{align}
\end{thm}

\proof ~For any initial vector $x_0\in R^n$, from Lemma 3 there has $x_k\in D_r(k\ge 0)$, thus $x_k-P_{N(A)}x_0\in N(A)^\bot$, therefore
\begin{align}
  x_{k+1}-P_{N(A)}x_0&=x_k-P_{N(A)}x_0+\frac{b_i-(a_i,x_k-P_{N(A)}x_0)}{\|a_i\|_2^2}a_i\nonumber\\
                     &=x_k-P_{N(A)}x_0+\frac{(Qb)_i-(a_i,x_k-P_{N(A)}x_0)}{\|a_i\|_2^2}a_i+\frac{b_i-(Qb)_i}{\|a_i\|_2^2}a_i.
\end{align}
Assume $x^\dagger$ is the generalized solution of $Ax=b$, hence
\begin{align}\label{error_expression}
  x_{k+1}-P_{N(A)}x_0-x^\dagger&=x_k-P_{N(A)}x_0-x^\dagger+\frac{(a_i,x^\dagger)-(a_i,x_k-P_{N(A)}x_0)}{\|a_i\|_2^2}a_i+\frac{\big{(}(I-Q)b\big{)}_i}{\|a_i\|_2^2}a_i\nonumber\\
                               &=x_k-P_{N(A)}x_0-x^\dagger-\frac{(a_i,x_k-P_{N(A)}x_0-x^\dagger)}{\|a_i\|_2^2}a_i+\frac{\big{(}(I-Q)b\big{)}_i}{\|a_i\|_2^2}a_i\nonumber\\
                               &=(I-\frac{a_ia_i^T}{\|a_i\|_2^2})(x_k-P_{N(A)}x_0-x^\dagger)+\frac{\big{(}(I-Q)b{)}_i}{\|a_i\|_2^2}a_i.
\end{align}
From Lemma \ref{orthogonal}, there holds
\begin{align}
 \big{(}(I-\frac{a_ia_i^T}{\|a_i\|_2^2})(x_k-P_{N(A)}x_0-x^\dagger),\frac{\big{(}(I-Q)b)_i}{\|a_i\|_2^2}a_i\big{)}=0,
\end{align}
consequently,
\begin{align}\label{eq.norm.1}
   \|x_{k+1}-P_{N(A)}x_0-x^\dagger\|_2^2=\|(I-\frac{a_ia_i^T}{\|a_i\|_2^2})(x_k-P_{N(A)}x_0-x^\dagger)\|_2^2+\|\frac{\big{(}(I-Q)b\big{)}_i}{\|a_i\|_2^2}a_i\|_2^2.
\end{align}
In addition,
\begin{align}\label{eq.norm.2}
    \|(I-\frac{a_ia_i^T}{\|a_i\|_2^2})(x_k-P_{N(A)}x_0-x^\dagger)\|_2^2=\|x_k-P_{N(A)}x_0-x^\dagger\|_2^2-\|\frac{a_i^T}{\|a_i\|_2}(x_k-P_{N(A)}x_0-x^\dagger)\|_2^2.
\end{align}
From \eqref{eq.norm.1} and \eqref{eq.norm.2}, there holds
\begin{align}\label{eq.norm}
   \|x_{k+1}-P_{N(A)}x_0-x^\dagger\|_2^2=&\|x_k-P_{N(A)}x_0-x^\dagger\|_2^2-\|\frac{a_i^T}{\|a_i\|_2}(x_k-P_{N(A)}x_0-x^\dagger)\|_2^2\nonumber\\
                                         &+\frac{1}{\|a_i\|_2^2}\big{(}(I-Q)b\big{)}_i^2.
\end{align}
Consequently, taking expectation on both sides about index $i$ yields
\begin{align}\label{eq.norm.3}
   &E\Big{[}\|x_{k+1}-P_{N(A)}x_0-x^\dagger\|_2^2\big{|}x_k\Big{]}\nonumber\\
   &=\|x_k-P_{N(A)}x_0-x^\dagger\|_2^2-\sum\limits_{i=1}^m\frac{\|a_i\|_2^2}{\|A\|_F^2}\|\frac{a_i^T}{\|a_i\|_2}(x_k-P_{N(A)}x_0-x^\dagger)\|_2^2+\sum\limits_{i=1}^m\frac{\|a_i\|_2^2}{\|A\|_F^2}\|\frac{1}{\|a_i\|_2^2}\big{(}(I-Q)b\big{)}_i^2\nonumber\\
   &=\|x_k-P_{N(A)}x_0-x^\dagger\|_2^2-\sum\limits_{i=1}^m\frac{\|a_i\|_2^2}{\|A\|_F^2}\frac{1}{\|a_i\|_2^2}\|a_i^T(x_k-P_{N(A)}x_0-x^\dagger)\|_2^2+\sum\limits_{i=1}^m\frac{1}{\|A\|_F^2}\big{(}(I-Q)b)_i^2\nonumber\\
   &=\|x_k-P_{N(A)}x_0-x^\dagger\|_2^2-\frac{1}{\|A\|_F^2}\|A(x_k-P_{N(A)}x_0-x^\dagger)\|_2^2+\frac{1}{\|A\|_F^2}\|(I-Q)b\|_2^2.
\end{align}
Since $x_k-P_{N(A)}x_0-x^\dagger\in N(A)^\bot$, there has
\begin{align}
  A^{\dagger}A(x_k-P_{N(A)}x_0-x^\dagger)=x_k-P_{N(A)}x_0-x^\dagger,
\end{align}
hence
\begin{align}
  \|x_k-P_{N(A)}x_0-x^\dagger\|_2^2=\| A^{\dagger}A(x_k-P_{N(A)}x_0-x^\dagger)\|_2^2\le \|A^\dagger\|_2^2\|A(x_k-P_{N(A)}x_0-x^\dagger)\|_2^2,
\end{align}
which yields
\begin{align}\label{eq.norm.4}
  \|A(x_k-P_{N(A)}x_0-x^\dagger)\|_2^2\ge\frac{1}{\|A^\dagger\|_2^2}\|x_k-P_{N(A)}x_0-x^\dagger\|_2^2.
\end{align}
Substitute \eqref{eq.norm.4} into \eqref{eq.norm.3}, there holds
\begin{align}\label{convergence_m}
  E\left[\|x_{k+1}-P_{N(A)}x_0-x^\dagger\|_2^2\right]&\le (1-\frac{1}{\|A\|_F^2\|A^\dagger\|_2^2})E\left[\|(x_k-P_{N(A)}x_0-x^\dagger)\|_2^2\right]+\frac{1}{\|A\|_F^2}\|(I-Q)b\|_2^2\nonumber\\
                                       &\le (1-\frac{1}{\|A\|_F^2\|A^\dagger\|_2^2})^k\|x_0-P_{N(A)}x_0-x^\dagger\|_2^2+\frac{1}{\|A\|_F^2}\|(I-Q)b\|_2^2.
\end{align}
Denotes $\kappa_A=\|A\|_F\|A\|_2$ as $\kappa(A)$, then \eqref{convergence_m} implies \eqref{thm.result.1}.\hfill{$\square$}

In \eqref{thm.result.1}, if $b\in R(A)$, it means $\|(I-Q)b\|_2\equiv 0$, then there holds
\begin{align*}
  E\left[\|x_{k+1}-P_{N(A)}x_0-x^\dagger\|_2^2\right]&\le (1-\frac{1}{\kappa^2(A)})^k\|x_0-P_{N(A)}x_0-x^\dagger\|_2^2.
\end{align*}
As $k\rightarrow \infty$, $1-\frac{1}{\kappa^2(A)}$ tends to zero, thus $\|x_{k+1}-P_{N(A)}x_0-x^\dagger\|_2^2\rightarrow 0$, i.e.
\begin{align*}
  x_{k+1}\rightarrow P_{N(A)}x_0+x^\dagger,
\end{align*}
which proves that the vector sequence $\{x_k\}_{k=1}^\infty$ generated from the randomized Kaczmarz method is convergent when $Ax=b$ is consistent and the limit is $P_{N(A)}x_0+x^\dagger$.

In addition, it is easy to see from theorem \ref{thm.5} that if $b\notin R(A)$, i.e. the system $Ax=b$ is inconsistent. The vector sequence generated by randomized Kaczmarz method is bounded rather than convergent(see \cite{Popa1995}).

In fact, Popa pointed out that the (classical) Kaczmarz algorithm with $x_0=0$ generates a sequence ${x_k}$ convergent to $x^\dagger$ if and only if the system $Ax=b$ is consistent.
The (classic) Kaczmrz and randomized Kaczmarz method are essentially the same in convergence behavior regardless of convergent rate.
Theorem \ref{thm.5} confirms that the Kaczmarz and randomized Kaczamrz methods for solving the insistent system do not converge to $x^\dagger$ even if $x_0=0$.

Unless the above illustrations,  error scheme \eqref{thm.result.1} can be used to analyze the perturbation problems when $\|(I-Q)b\|_2$ is replaced with $\|b^\delta-b\|_2$, where $b^\delta$ is the noisy right-hand side of the consistent system $Ax=b$. Kaczmarz like methods for solving ill-posed problems have semi-convergent behavior (see \cite{Elfing_2014}), however, Theorem \ref{thm.5} show us an upper bound of $\|x_k-P_{N(x)x_0}-x^\dagger\|_2$.

During the performance of randomized Kaczmarz method \eqref{randomized_kaczmarz}, the probability to choose the $i$th equation is $\frac{\|a_i\|_2^2}{\|A\|_F^2}$, but in fact, if the normal vector of every equation is canonical, i.e. $\|a_i\|_2=1, i=1,2,\cdots,m$, then the probability in randomized Kaczmarz method is changed to $\frac{1}{m}$. For general linear system, we can acquire it by normalizing the coefficient matrix $A$, so the original linear system $Ax=b$ is equivalent to $DAx=Db$, where $D=\text{diag}(1/\|a_1\|_2,1/\|a_2\|_2,\cdots,1/\|a_m\|_2)$. From Theorem \ref{thm.randomized.Kaczmarz}, there holds the next corollary for the normalized linear system.
\begin{cor}\label{normal_system.cor}\normalfont
  For linear system $Ax=b$, $D=\text{diag}(1/\|a_1\|_2,1/\|a_2\|_2,\cdots,1/\|a_m\|_2)$, $P_{N(A)}: R^n\rightarrow N(A)$ and $Q: R^m\rightarrow R(A)$ are orthogonal projection operators, $\{x_k\}_{k=1}^\infty$ is a vector sequence generated by randomized Kaczmarz method, then there holds\\
  (1) For inconsistent data,
      \begin{align}\label{normal_result.1}
          E\big{[}\|x_{k+1}-P_{N(A)}x_0-x^\dagger\|_2^2\big{]}\le&(1-\frac{1}{m\cdot\max\limits_{i=1,\cdots,m}\|a_i\|_2^2\|A^\dagger\|_2^2})^k\|x_0-P_{N(A)}x_0-x^\dagger\|_2^2\nonumber\\
                                                                 &+\frac{1}{\|A\|_F^2}\|(I-Q)Db\|_2^2.
      \end{align}
  (2) For consistent data,
      \begin{align}\label{normal_result.2}
          E\big{[}\|x_{k+1}-P_{N(A)}x_0-x^\dagger\|_2^2\big{]}\le (1-\frac{1}{m\cdot\max\limits_{i=1,\cdots,m}\|a_i\|_2^2\|A^\dagger\|_2^2})^k\|x_0-P_{N(A)}x_0-x^\dagger\|_2^2.
      \end{align}
\end{cor}
\proof ~~For linear system $DAx=Db$, from the theorem \ref{thm.5}, hence
\begin{align*}
   E\big{[}\|x_{k+1}-P_{N(A)}x_0-x^\dagger\|_2^2\big{]}\le (1-\frac{1}{\|DA\|_F^2\|(DA)^\dagger\|_2^2})^k\|x_0-P_{N(A)}x_0-x^\dagger\|_2^2+\frac{1}{\|A\|_F^2}\|(I-Q)b\|_2^2.
\end{align*}
Note that $\|DA\|_F^2=m, \|D^{-1}\|_2^2=\max(\|a_i\|_2^2)$, therefore
\begin{align*}
    \|DA\|_F^2\|(DA)^\dagger\|_2^2\le m\max\limits_i\|a_i\|_2^2\|A^\dagger\|_2^2.
\end{align*}
Consequently,
\begin{align*}
   E\big{[}\|x_{k+1}-P_{N(A)}x_0-x^\dagger\|_2^2\big{]}\le (1-\frac{1}{m\max\limits_{i=1,\cdots,m}\|a_i\|_2^2\|A^\dagger\|_2^2})^k\|x_0-P_{N(A)}x_0-x^\dagger\|_2^2+\frac{1}{\|A\|_F^2}\|(I-Q)Db\|_2^2.
\end{align*}

Furthermore, from $\|A\|_2^2=\sum\limits_{i=1}^m\|a_i\|_2^2\le m\max\|a_i\|_2^2$, there holds
\begin{align}
  (1-\frac{1}{\|A\|_F^2\|A^\dagger\|_2^2})\le(1-\frac{1}{m\max\limits_{i=1,\cdots,m}\|a_i\|_2^2\|A^\dagger\|_2^2}).
\end{align}
Therefore, the result of the Corollary \ref{normal_system.cor} for normal system is slightly weaker than the result of Theorem \ref{thm.5} for general linear system.
In fact, the result can be understood easily, in order to keep the consistency with $Ax=b$, we hope to express the result with $A$ or some parts of it rather than $DA$, so it is unavoidable to enlarge appropriately the item $\|A\|_F\|A^\dagger\|_2$ in \eqref{thm.result.1}.

For the consistent system, i.e. $(I-Q)b=0$,  the inequality \eqref{normal_result.1} can be simplified to \eqref{normal_result.2}.\hfill{$\square$}

Next, we consider the classical Kaczmarz method, there holds the following error estimate.
\begin{thm}\label{thm.6}
  \normalfont
  If $Ax=b$ is consistent, denote $N_{k+1}=N(a_{k+1})$. $P_{N(A)}: R^n\mapsto N(A)^\bot$ and $P_{k+1}:N(A)^\bot\mapsto N_{k+1}^\bot$ are orthogonal projection operators. For any
initial vector $x_0\in R^n$, the vector consequence $\{x_k\}_{k=1}^m$ generated from Kaczmarz method \eqref{kaczmarz.1} satisfies
\begin{align}\label{thm.result.2}
   \|x_{k+1}-P_{N(A)}x_0-x^\dagger\|_2^2\le (1-\frac{1}{\|a_{k+1}\|_2^2\|(a_{k+1}^TP_{k+1})^\dagger\|_2^2})\|x_k-P_{N(A)}x_0-x^\dagger\|_2^2
\end{align}
and $\|(a_{k+1}^TP_{k+1})^\dagger\|_2\ge \|A^\dagger\|_2$.
\end{thm}
\proof ~~In \eqref{eq.norm}, let $i=k+1$. Since $Ax=b$ is consistent, thus $Q=I$, and hence
\begin{align}\label{eq.norm.5}
   \|x_{k+1}-P_{N(A)}x_0-x^\dagger\|_2^2&=\|x_k-P_{N(A)}x_0-x^\dagger\|_2^2-\|\frac{a_{k+1}^T}{\|a_{k+1}\|_2}(x_k-P_{N(A)}x_0-x^\dagger)\|_2^2\nonumber\\
                                        &\le\|x_k-P_{N(A)}x_0-x^\dagger\|_2^2-\|\frac{a_{k+1}^T}{\|a_{k+1}\|_2}P_{k+1}(x_k-P_{N(A)}x_0-x^\dagger)\|_2^2.
\end{align}
for $N(A)\subset N_{k+1}$, thus $N_{k+1}^\bot\subset N(A)^\bot$, from the condition $P_{k+1}:N(A)^\bot\mapsto N_{k+1}^\bot$, let $(a_{k+1}^TP_{k+1})^\dagger$ be the generalized inverse of $a_{k+1}^TP_{k+1}$, there holds
\begin{align}\label{eq.norm.6}
  \|a_{k+1}^TP_{k+1}(x_k-P_{N(A)}x_0-x^\dagger)\|_2^2\ge \frac{1}{\|(a_{k+1}^TP_{k+1})^\dagger\|_2^2}\|x_k-P_{N(A)}x_0-x^\dagger\|_2^2,
\end{align}
substitute \eqref{eq.norm.5} into \eqref{eq.norm.6}, then there holds
\begin{align}
   \|x_{k+1}-P_{N(A)}x_0-x^\dagger\|_2^2\le(1-\frac{1}{\|a_{k+1}\|_2^2\|(a_{k+1}^TP_{k+1})^\dagger\|_2^2})\|x_k-P_{N(A)}x_0-x^\dagger\|_2^2,
\end{align}
the inequality \eqref{thm.result.2} is proved.

In addition,
\begin{align*}
  \|(a_{k+1}^TP_{k+1})^\dagger\|_2&=\sup\limits_{y\in R(a_{k+1})}\|\frac{(a_{k+1}^TP_{k+1})^\dagger y\|_2}{\|y\|_2}=\sup\limits_{x\in N_A^\bot}\frac{\|(a_{k+1}^TP_{k+1})^\dagger a_{k+1}^TP_{k+1}x\|_2}{\|a_{k+1}^TP_{k+1}x\|_2}\nonumber\\
  &=\sup\limits_{x\in N_A^\bot}\frac{\|x\|_2}{\|a_{k+1}^TP_{k+1}x\|_2}
  =\sup\limits_{x\in N_A^\bot}\frac{\|x\|_2}{\|a_{k+1}^Tx\|_2}\ge \sup\limits_{x\in N_A^\bot}\frac{\|x\|_2}{\|Ax\|_2}=\|A^\dagger\|_2,
\end{align*}

It is obvious that the above proof about the error estimate for the vector sequence generated by classical Kaczmarz method is restricted in one recycle period, i.e. $k=1,2,\cdots,m$, in other word, the index of error is in accordance with the number of equation. For general case, there holds the following result.

\begin{cor}\label{corollary.7}
  \normalfont
  If $Ax=b$ is consistent, denote $N_{k+1}=\ker(a_{k+1})=N(a_{k+1})$, $P_{N(A)}: R^n\mapsto N(A)^\bot$ and $P_{k+1}:N(A)^\bot\mapsto N_{k+1}^\bot$ are orthogonal projection operators, for any initial vector $x_0\in R^n$, the vector consequence $\{x_k\}_{k=1}^\infty$ generated from Kaczmarz method \eqref{kaczmarz.1} satisfies
\begin{align}\label{cor.result.2}
   \|x_{k+1}-P_{N(A)}x_0-x^\dagger\|_2^2\le (1-\frac{1}{\max\limits_{i=1,2,\cdots,m}\|a_i\|_2^2\|(a_i^TP_i)^\dagger\|_2^2})^{k+1}\|x_0-P_{N(A)}x_0-x^\dagger\|_2^2.
\end{align}
\end{cor}

From Corollary \ref{corollary.7}, as $k\rightarrow \infty$, there holds $\|x_{k+1}-P_{N(A)}x_0-x^\dagger\|_2\rightarrow 0$ at once, i.e., $x_{k+1}\rightarrow P_{N(A)}x_0+x^\dagger$ which is in accordance with the conclusion of Kunio Tanabe in \cite{Tanabe_1971}.

\section{Block randomized Kaczmarz method}

In this part, we will consider the block randomized Kaczmrz method. Assume the linear system $Ax=b$, where $A\in R^{m\times n}, b\in R^m$, $S=\{1,2,\cdots,m\}$ is the identifier set of the equations of the linear system, divide the set $S=S_1\cup S_2\cup \cdots \cup S_r$, where $S_i\cap S_j=\Phi(i,j=1,2,\cdots,r)$, actually, $\{S_i\}_{i=1}^r$ is a classification set of $S$.  the number of the elements in $S_i$ is denoted by $\#S_i$.

Based on the general framework, i.e. without fixing classification set, The algorithm of the Block randomized Kaczmarz method can be provided in the following.

\begin{algorithm}
\caption{Block randomized Kaczmarz method}\label{block.randomized.Kaczmarz} 
\begin{algorithmic}[1]
\State Given $A,b,x_0$, $S_1,S_2,\cdots,S_r$, $N$.
\State $k=1$, let $S_{sel}=S_1$.
\State perform
   \begin{align}\label{iteration.1}
      x_k=x_{k-1}+\frac{b_i-(a_i,x_{k-1})}{\|a_i\|_2^2}a_i,
   \end{align}
where the identifier $i$ is selected by the probability $P_r=\frac{1}{\#S_{sel}}$ in subset $S_{sel}$.
\State $k=k+1$, if $k<N$, let $l=k\bmod r$, if $l\neq 0$, let $S_{sel}=S_l$ and if $l=0$, let $S_{sel}=S_r$, then select one number $i$ in $S_{sel}$ by the probability $\frac{1}{\#S_l}$ and then go to step 3; otherwise, terminate the iteration and let $x_k$ be the numerical solution.
\end{algorithmic}
\end{algorithm}

For block randomized Kaczmarz method, there holds the next theorem.
\begin{thm}\label{theorem.blocker.randomized.kaczamrz}\normalfont
  Assume identifer $i$ is selected in $S_c$, where $c$ is a certain value in set $\{1,2,\cdots,r\}$, $P_i: N(A)^\bot\mapsto N(a_i)^\perp$, then there holds for the vector sequence $\{x_k\}_{k=1}^\infty$ generated by Algorithm \ref{block.randomized.Kaczmarz}:

  (i) For consistent system
      \begin{align}
          E\Big{[}\|x_{k+1}-P_{N(A)}x_0-x^\dagger\|_2^2\big{|}x_k\Big{]}\le (1-\min\limits_{i\in S_c}\frac{1}{\|a_i\|_2^2\|(a_i^TP_i)^\dagger\|_2^2})\|x_k-P_{N(A)}x_0-x^\dagger\|_2^2.
      \end{align}

  (ii) For inconsistent system
      \begin{align}
          E\big{[}\|x_{k+1}-P_{N(A)}x_0-x^\dagger\|_2^2\big{]}\le (1-\frac{1}{m\max\limits_{i=1,\cdots,m}\|a_i\|_2^2\|A^\dagger\|_2^2})^k\|x_0-P_{N(A)}x_0-x^\dagger\|_2^2+\frac{m}{\#S_{sel}}\frac{\delta^2}{\lambda_{\min}^2(A)}.
      \end{align}
\end{thm}
\proof~~From \eqref{eq.norm}
\begin{align*}
   \|x_{k+1}-P_{N(A)}x_0-x^\dagger\|_2^2=\|x_k-P_{N(A)}x_0-x^\dagger\|_2^2-\|\frac{a_i^T}{\|a_i\|_2}(x_k-P_{N(A)}x_0-x^\dagger)\|_2^2+\frac{1}{\|a_i\|_2^2}\big{(}(I-Q)b)_i^2.
\end{align*}
Taking expectation on both side about $i$ in $S_c$ yields
\begin{align*}
   &E\Big{[}\|x_{k+1}-P_{N(A)}x_0-x^\dagger\|_2^2\big{|}x_k\Big{]}\nonumber\\
   &=\|x_k-P_{N(A)}x_0-x^\dagger\|_2^2-\frac{1}{\#S_c}\sum\limits_{i=1}^{\#S_c}\|\frac{a_i^T}{\|a_i\|_2}(x_k-P_{N(A)}x_0-x^\dagger)\|_2^2+\frac{1}{\#S_c}\sum\limits_{i=1}^{\#S_c}\frac{1}{\|a_i\|_2^2}\big{(}(I-Q)b)_i^2\nonumber\\
   &\le\|x_k-P_{N(A)}x_0-x^\dagger\|_2^2-\min\limits_{i\in S_c}\frac{1}{\|a_i\|_2^2}\|a_i^T(x_k-P_{N(A)}x_0-x^\dagger)\|_2^2+\frac{1}{\#S_c}\frac{1}{\min\limits_{i\in S_c}\|a_i\|_2^2}\|(I-Q)b\|_2^2\nonumber\\
   &\le\|x_k-P_{N(A)}x_0-x^\dagger\|_2^2-\min\limits_{i\in S_c}\frac{1}{\|a_i\|_2^2}\|a_i^TP_i(x_k-P_{N(A)}x_0-x^\dagger)\|_2^2+\frac{1}{\#S_c}\frac{1}{\min\limits_{i\in S_c}\|a_i\|_2^2}\|(I-Q)b\|_2^2\nonumber\\
   &\le(1-\min\limits_{i\in S_c}\frac{1}{\|a_i\|_2^2\|(a_i^TP_i)^\dagger\|_2^2})\|x_k-P_{N(A)}x_0-x^\dagger\|_2^2+\frac{1}{\#S_c}\frac{1}{\min\limits_{i\in S_c}\|a_i\|_2^2}\|(I-Q)b\|_2^2.
\end{align*}
The inequality \eqref{block.inconsistent} is proved. Taking $Qb=b$ will obtain the inequality for the consistent system.

From Theorem \ref{theorem.blocker.randomized.kaczamrz}, the following corollary is obvious.
\begin{cor}\label{cor.block.randomized.Kaczmarz}\normalfont

Assume $Q: R^m\rightarrow R(A)$ and  $P_i:N(A)^\bot\mapsto N(a_i)^\bot$ are orthogonal projection operators,  then there hold for the vector sequence $\{x_k\}_{k=1}^\infty$ generated by Algorithm \ref{block.randomized.Kaczmarz}:

  (i) For consistent system
      \begin{align}\label{block.consistent}
          E[\|x_{k+1}-P_{N(A)}x_0-x^\dagger\|_2^2]\le(1-\min\limits_{i\in S}\frac{1}{\|a_i\|_2^2\|(a_i^TP_i)^\dagger\|_2^2})^{k+1}\|x_0-P_{N(A)}x_0-x^\dagger\|_2^2.
      \end{align}

  (ii) For inconsistent system
      \begin{align}\label{block.inconsistent}
          E[\|x_{k+1}-P_{N(A)}x_0-x^\dagger\|_2^2]\le&(1-\min\limits_{i\in S}\frac{1}{\|a_i\|_2^2\|(a_i^TP_i)^\dagger\|_2^2})^{k+1}\|x_0-P_{N(A)}x_0-x^\dagger\|_2^2\nonumber\\
                      &+\frac{\max\limits_{i\in S}\|a_i\|_2^2}{\min\limits_{i\in S}\|a_i\|_2^2}\frac{\max\limits_{i\in S}\|(a_i^TP_i)^\dagger\|_2^2}{\min\limits_{c\in\{1,2,\cdots,r\}}\#S_c}\|(I-Q)b\|_2^2.
      \end{align}
\end{cor}

 when the projection is performed randomly in classification sets $S_1,S_2,\cdots,S_r$ from $x_k$, the behavior of error for Kaczmarz method are shown in Theorem \ref{theorem.blocker.randomized.kaczamrz} and the corollary \ref{cor.block.randomized.Kaczmarz}.

\section{Numerical experiments}

In this section, we will illustrate these error estimate results appeared in the above sections by several classical problems, i.e. phillips, gravity and shaw. They are Fredholm integral equations of the first kind, phillips problem is mildly ill-posed and the others are severely ill-posed. If discretized them with dimension $m=n=1000$, the condition numbers of them are $2.6415e+10$, $9.8894e+19$ and $7.1301e+20$, respectively. the codes of discretized problems are taken from Matlab package Regutools\cite{Hansen1994}\footnote{Available from www.imm.dtu.dk/~pcha/Regutools/}.

In the test programs, discretized dimension is fixed to $m=n=1000$. We mainly compare the results of kaczmarz method and randomized Kaczmarz method. In perturbed case, the right-hand $b^\delta$ is generated from accurate term $b$, i.e.,
\begin{align*}
   b_i^\delta=b_i+\delta\max\limits_{i}(|b_i|),\qquad i=1,2,\cdots,n.
\end{align*}

\begin{figure}[ht]
  \centering
  \begin{minipage}[ht]{.25\linewidth}
      \includegraphics[width=120pt]{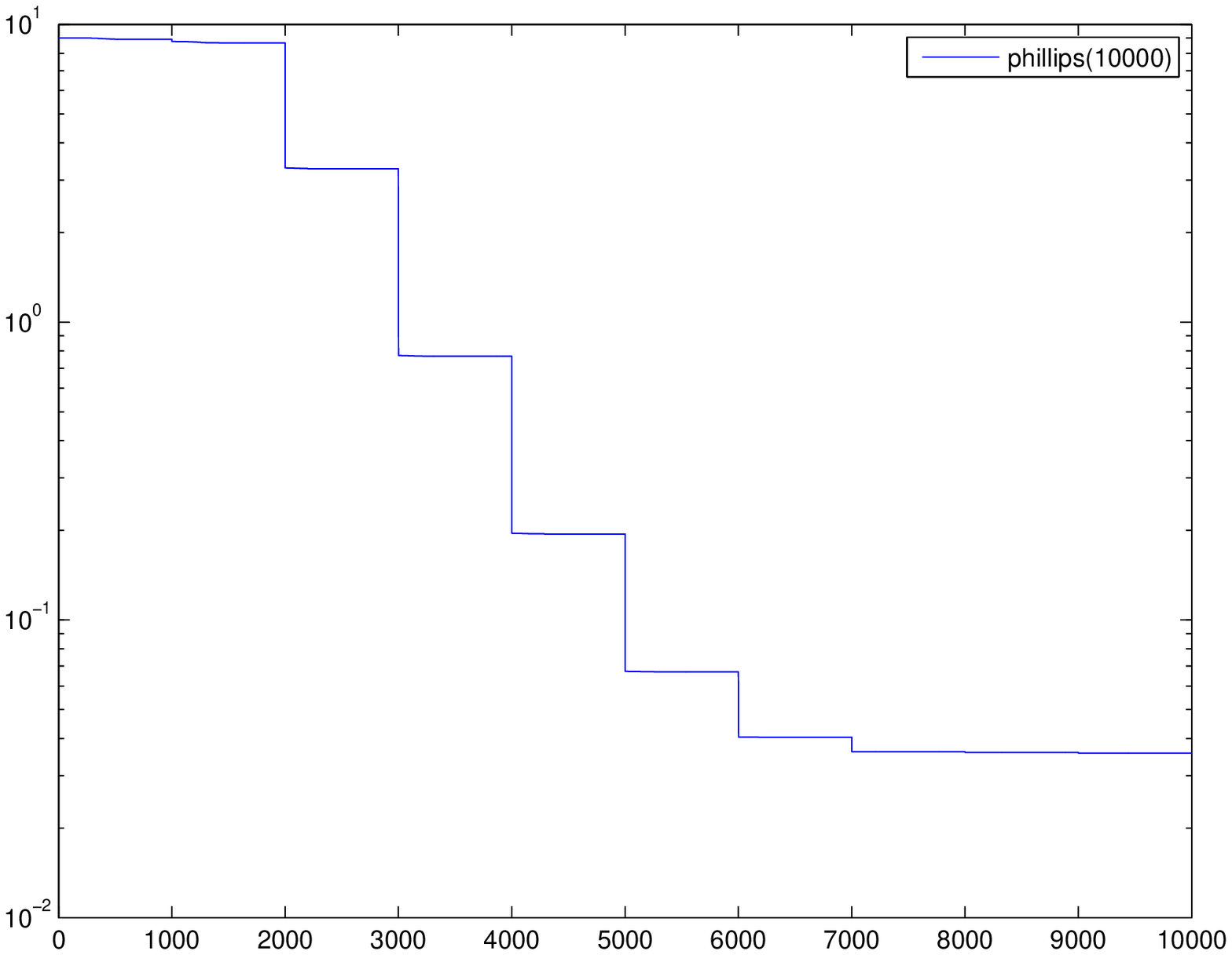}
       \caption{\footnotesize Kaczmarz method, K=10000,$\delta=0$}\label{kacphi(10000d0)}
  \end{minipage}
  \begin{minipage}[h]{.25\linewidth}
     \includegraphics[width=120pt]{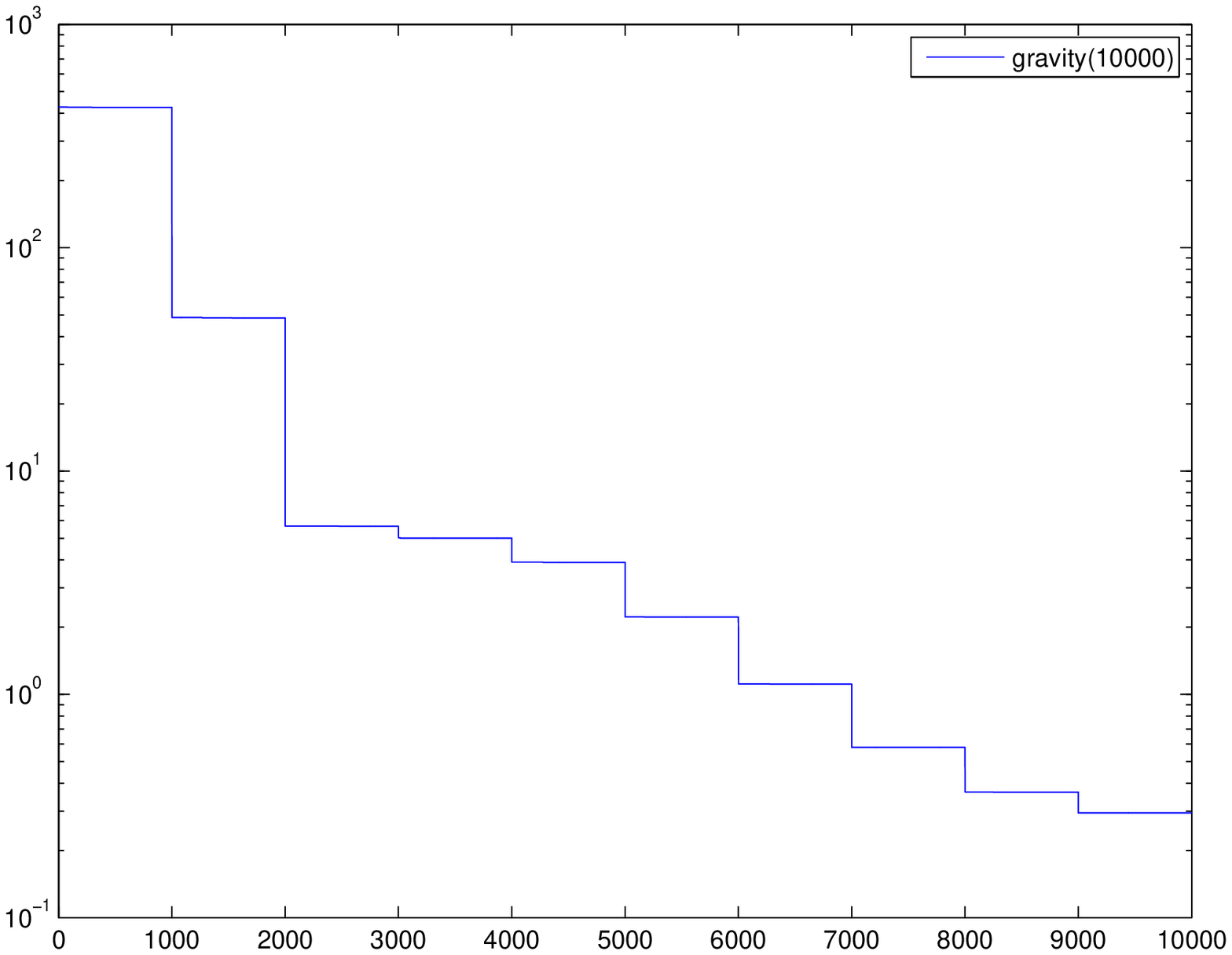}
     \caption{\footnotesize Kaczmarz method, K=10000,$\delta=0$}\label{kacgra(10000d0)}
  \end{minipage}
  \begin{minipage}[h]{.25\linewidth}
     \includegraphics[width=120pt]{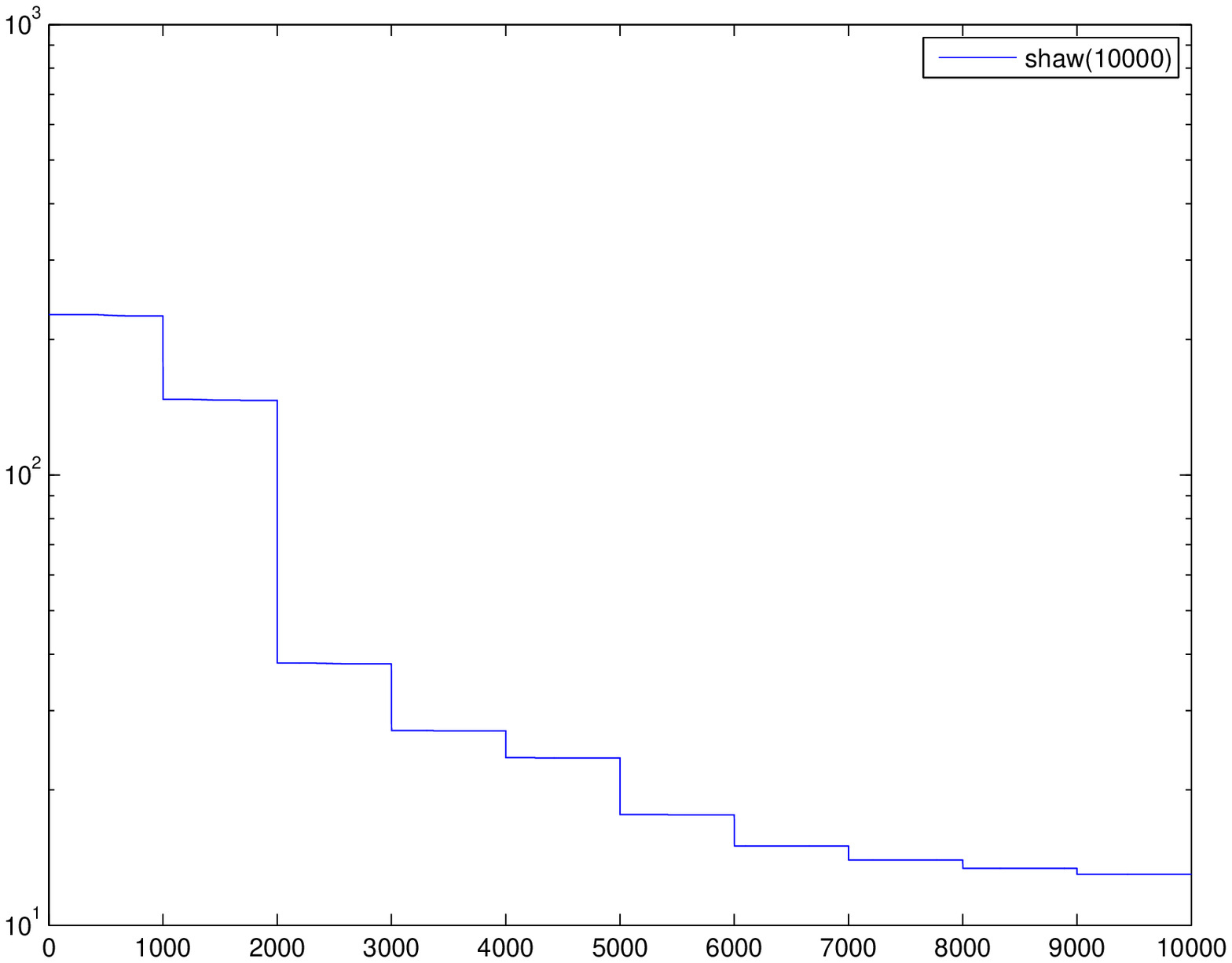}
         \caption{\footnotesize Kaczmarz method, K=10000,$\delta=0$}\label{kacsha(10000d0)}
 \end{minipage}
\end{figure}

Figure \ref{kacphi(10000d0)}$\sim$\ref{kacsha(10000d0)} are the figures of phillips, gravity and shaw solved by Kaczmarz as $K=10000, \delta=0$ and the initial vector $x_0=0$. From the behavior of the error estimate about the iterative step, it is easy to find that there are some 'ladders' between two recycles, in other words, the iterative improvement of Kaczmarz method is slight within one recycle until at the beginning of the next recycle. In fact, the behavior due to the high correlation between two adjacent equations.


\begin{figure}[ht]
  \centering
  \begin{minipage}[hb]{.25\linewidth}
      \includegraphics[width=120pt]{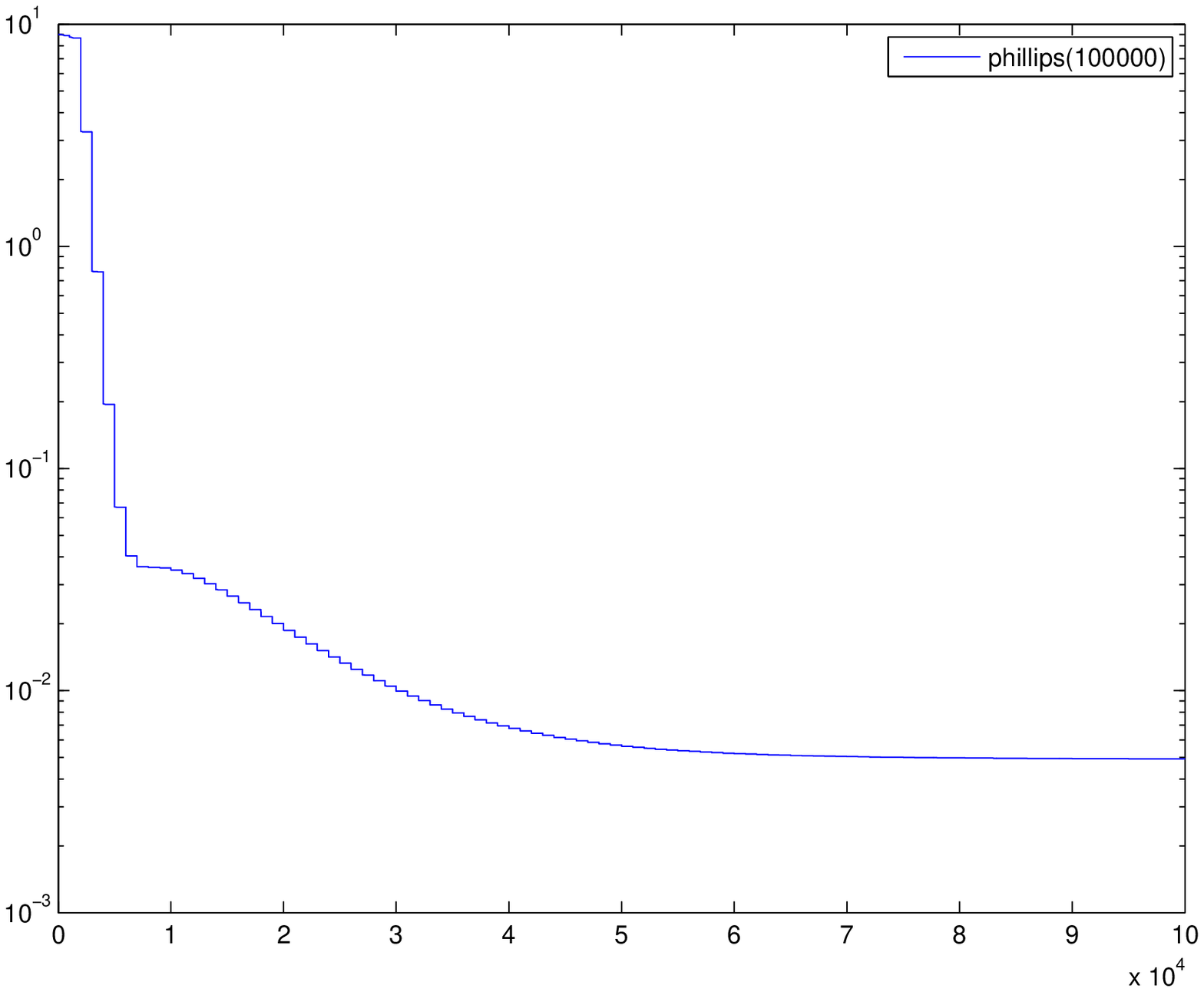}
       \caption{\scriptsize Kaczmarz method, K=100000,$\delta=0$}\label{kacphi(100000d0)}
  \end{minipage}
  \begin{minipage}[h]{.25\linewidth}
     \includegraphics[width=120pt]{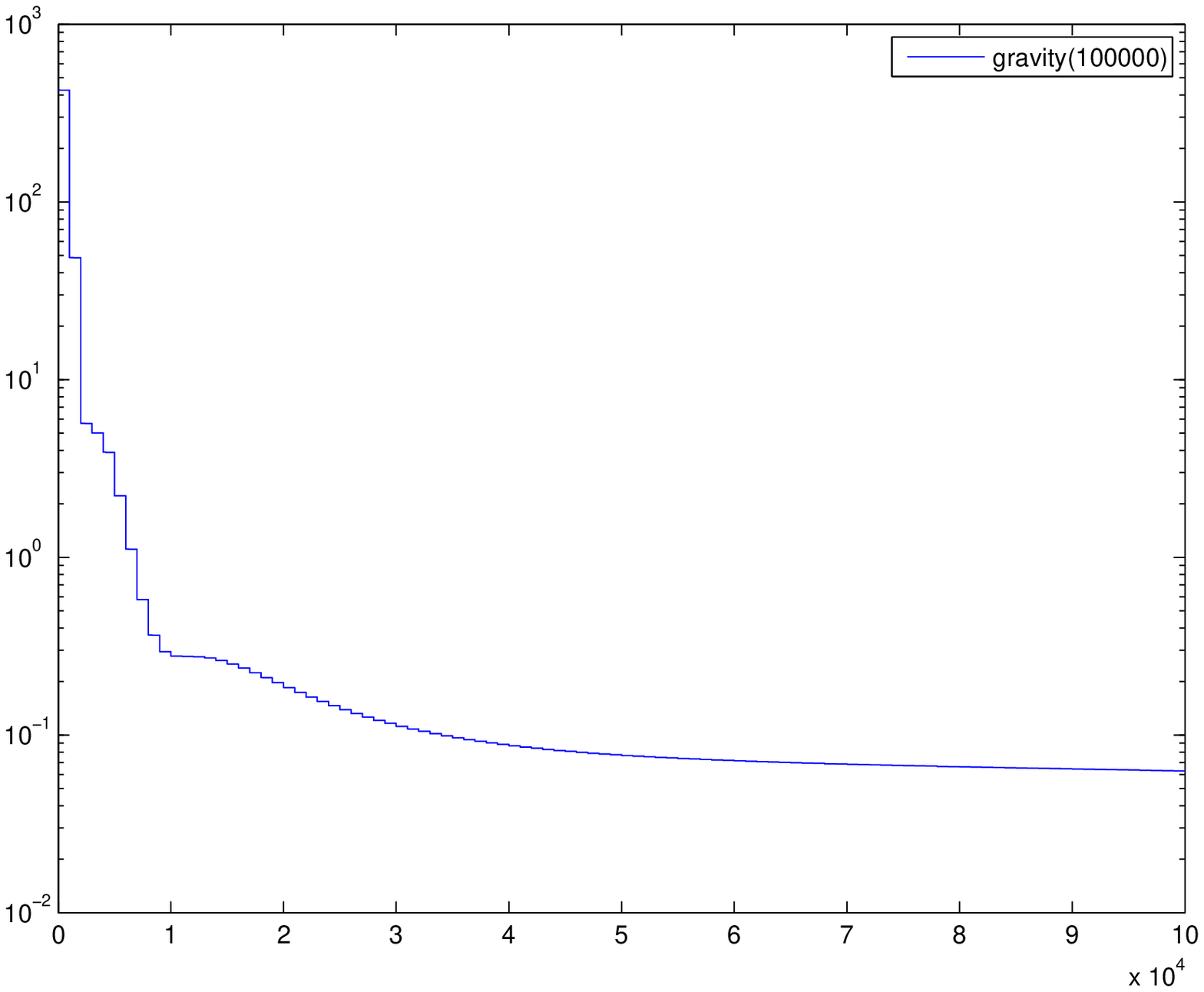}
     \caption{\scriptsize Kaczmarz method, K=100000,$\delta=0$}\label{kacgra(100000d0)}
  \end{minipage}
  \begin{minipage}[h]{.25\linewidth}
     \includegraphics[width=120pt]{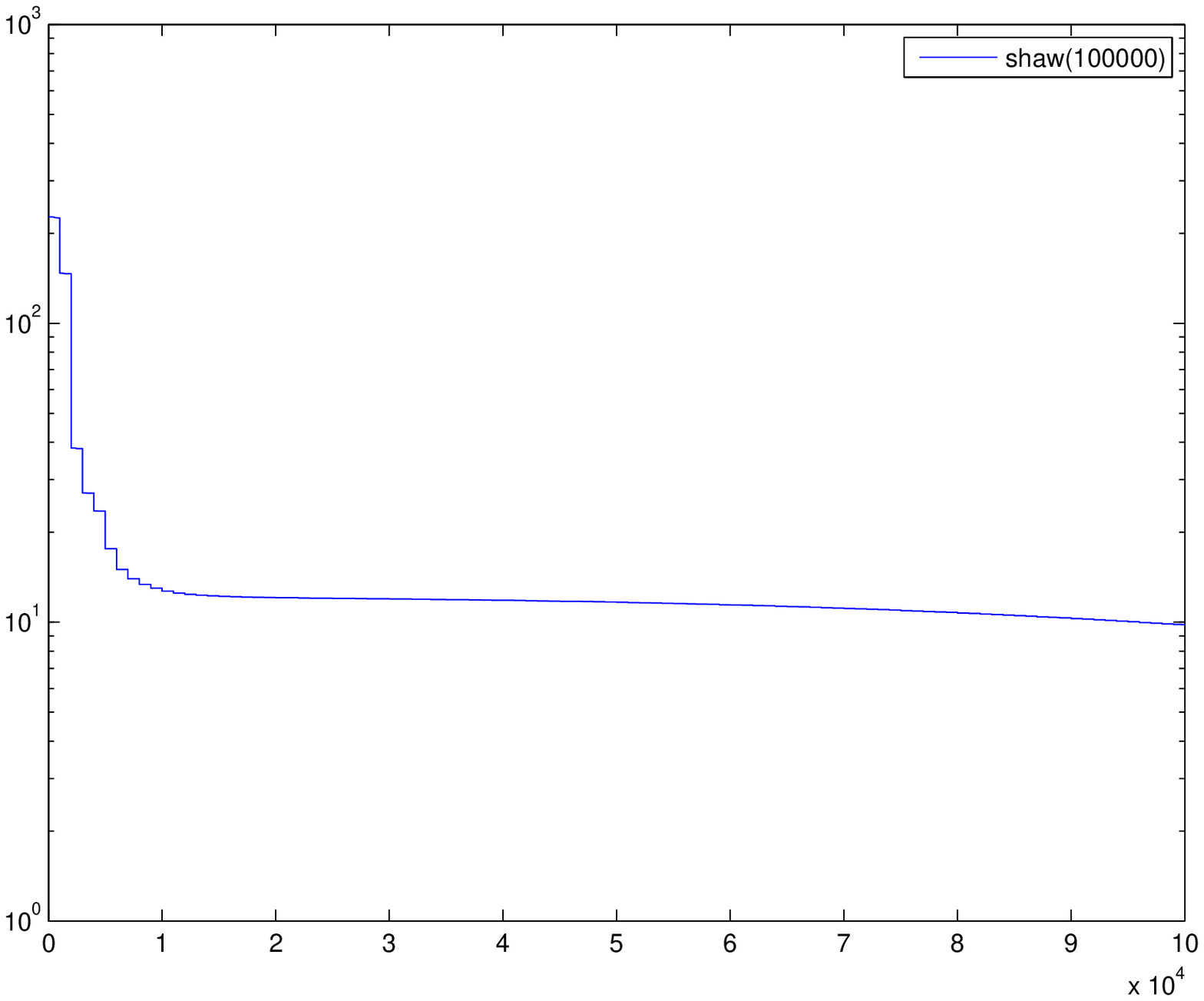}
         \caption{\scriptsize Kaczmarz method, K=100000,$\delta=0$}\label{kacsha(100000d0)}
 \end{minipage}
\end{figure}

In addition, we also illustrate the behavior from the following deduction. Let $x_k$ be the current iterative solution, we projected $x_k$ to the $(k+1)$th and the $(k+2)$th equations successively, i.e.
\begin{align*}
  &x_{k+1}=x_k+\frac{b_{k+1}-(a_{k+1},x_k)}{\|a_{k+1}\|_2^2}a_{k+1},\\
  &x_{k+2}=x_{k+1}+\frac{b_{k+2}-(a_{k+2},x_k)}{\|a_{k+2}\|_2^2}a_{k+2}.
\end{align*}
Let $x^*$ be an any solution of linear system $Ax=b$, and denote the error $e_{k+1}=x_{k+1}-x^*$, therefore,
\begin{align*}
  e_{k+1}-e_{k+2}&=\frac{a_{k+2}a_{k+2}^T}{\|a_{k+2}\|_2^2}(I-\frac{a_{k+1}a_{k+1}^T}{\|a_{k+1}\|_2^2})e_k\nonumber\\
  &=\frac{a_{k+2}}{\|a_{k+2}\|_2}(\frac{a_{k+2}^T}{\|a_{k+2}\|_2}-\frac{a_{k+1}^T\cos(\widehat{a_{k+1},a_{k+2}},)}{\|a_{k+1}\|_2})e_k,
\end{align*}
so,
\begin{align}\label{eq.no.1}
  \bigl{|}\|e_{k+1}\|-\|e_{k+2}\|\bigr{|}&\le \|e_{k+1}-e_{k+2}\|\le\bigl{\|}\frac{a_{k+2}^T}{\|a_{k+2}\|_2}-\frac{a_{k+1}^T\cos(\widehat{a_{k+1},a_{k+2}})}{\|a_{k+1}\|_2}\bigr{\|}\|e_k\|\nonumber\\
                           &\le\bigl{(\|}\frac{a_{k+2}^T}{\|a_{k+2}\|_2}\|_2^2+\|\frac{a_{k+1}^T}{\|a_{k+1}\|_2}\|_2^2
                           -2\cos(\widehat{a_{k+1},a_{k+2}})\langle\frac{a_{k+2}^T}{\|a_{k+2}\|_2},\frac{a_{k+1}^T}{\|a_{k+1}\|_2}\rangle\bigr{)}^\frac{1}{2}\|e_k\|\nonumber\\
                           &=\bigl{(}2-2\cos^2(\widehat{a_{k+1},a_{k+2}})\bigr{)}^{\frac{1}{2}}\|e_k\|.
\end{align}
From \eqref{eq.no.1}, if the intersection angle between two adjacent projection equations is small, the improvement of the current iteration $x_{k+1}$ is slight than the last iteration $x_k$, Because the intersection angle between the normal vectors of the adjacent equations are very small for these test problems, the behavior of their errors are as we can see in Figure \ref{kacphi(10000d0)}$\sim$\ref{kacsha(10000d0)}. Meanwhile,  these errors don't converge to zero, i.e. the numerical solutions don't converge to the original solutions(which can also be seen from Figure  \ref{kacphi(100000d0)}$\sim$\ref{kacsha(100000d0)}). In fact, from Theorem \ref{thm.1} and Theorem \ref{thm.6}, the numerical solutions of Kaczmarz method converges to their corresponding Moore-Penrose generalized solutions.

Let $K=100000$ and $\delta=0$, we can observe the subsequent behavior of these test problems solved by Kaczmarz method from Figure \ref{kacphi(100000d0)}$\sim$\ref{kacsha(100000d0)}, the error changes of these problems tend to be stable.

\begin{figure}[htb]
  \centering
  \begin{minipage}[h]{.25\linewidth}
      \includegraphics[width=130pt]{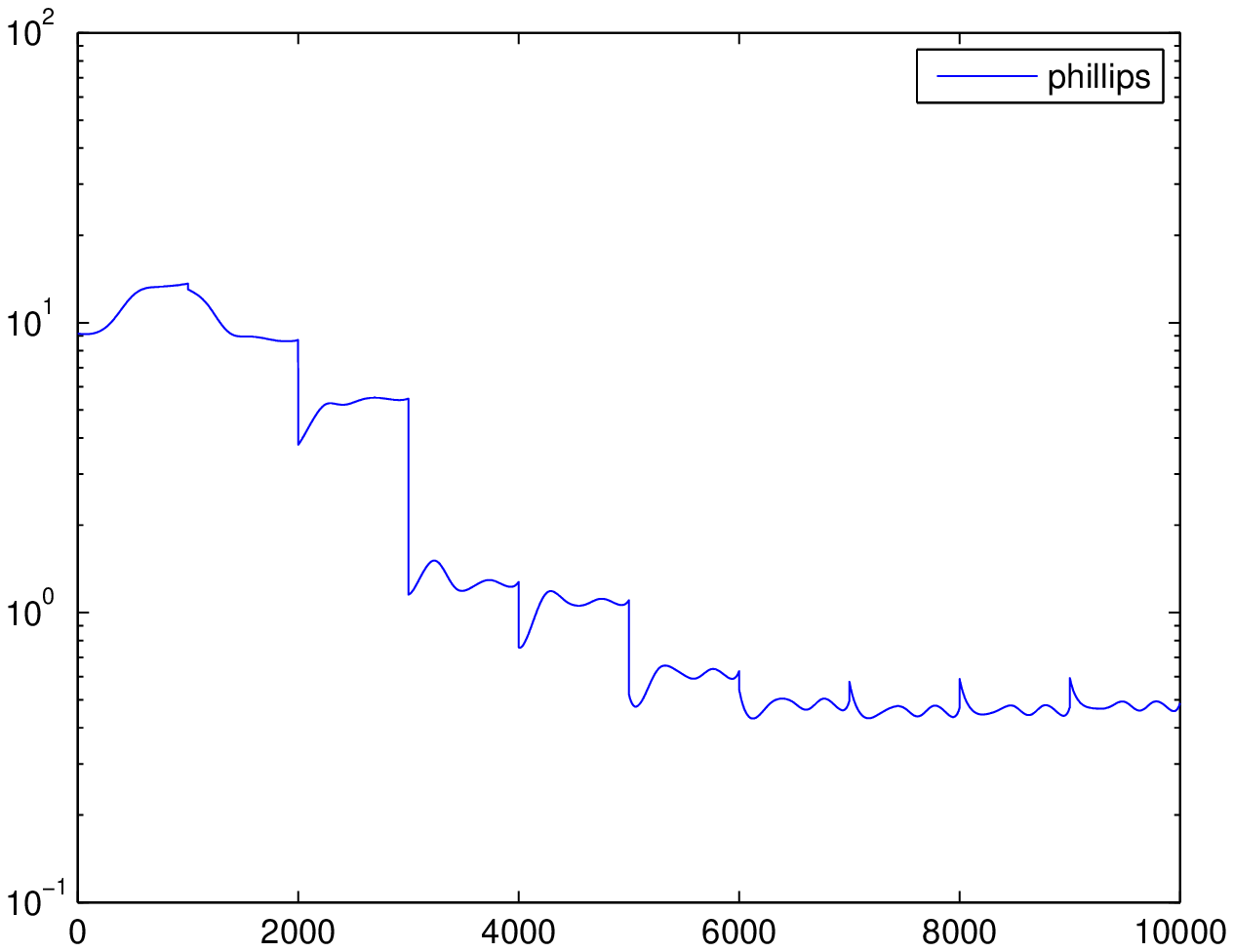}
       \caption{\footnotesize phillips(10000),\protect\\$\delta=0.1$}\label{kacphi(10000d0.1)}
  \end{minipage}
  \begin{minipage}[h]{.25\linewidth}
     \includegraphics[width=130pt]{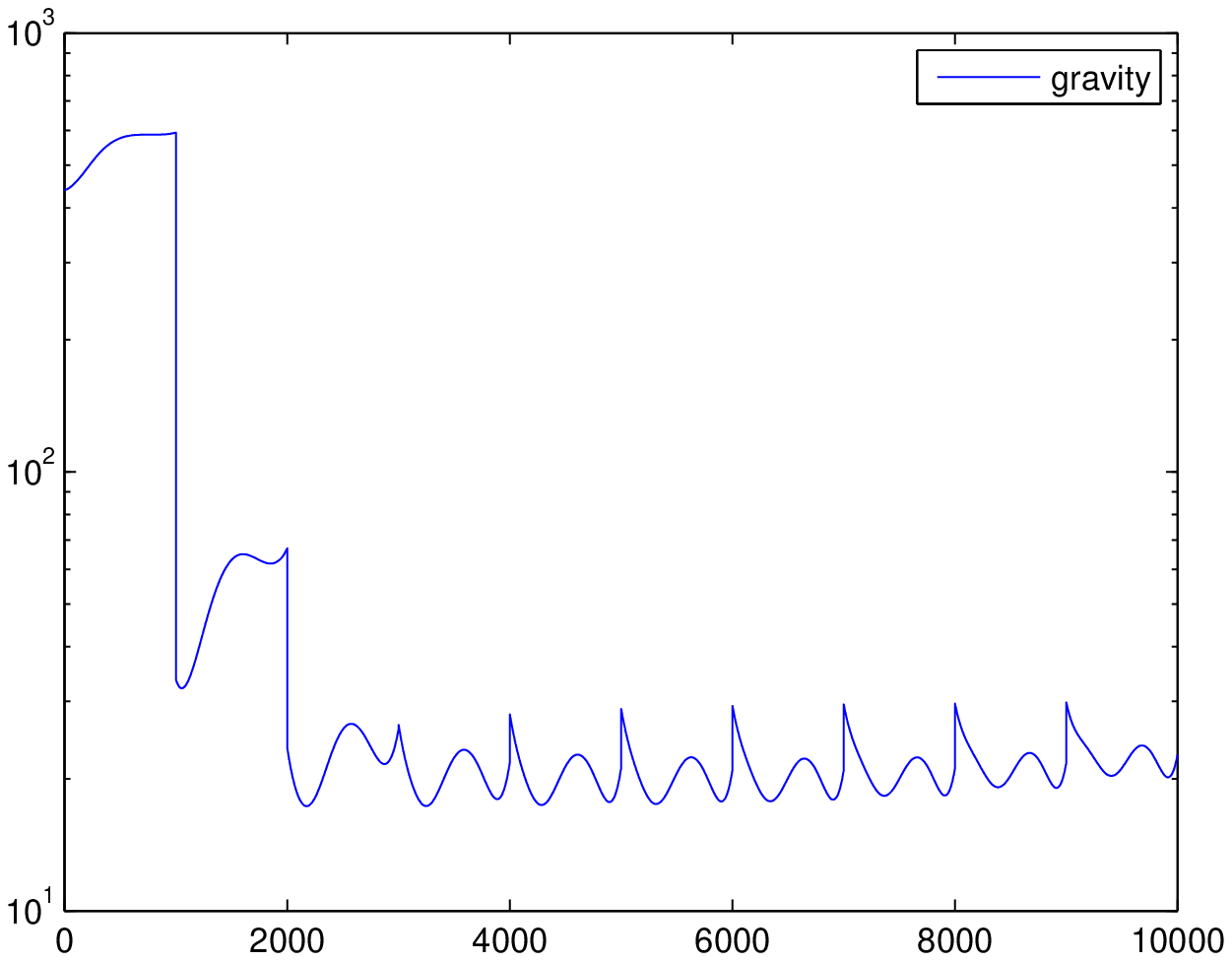}
     \caption{\footnotesize gravity(10000),\protect\\$\delta=0.1$}\label{kacgravity(10000d0.1)}
  \end{minipage}
  \begin{minipage}[h]{.25\linewidth}
     \includegraphics[width=130pt]{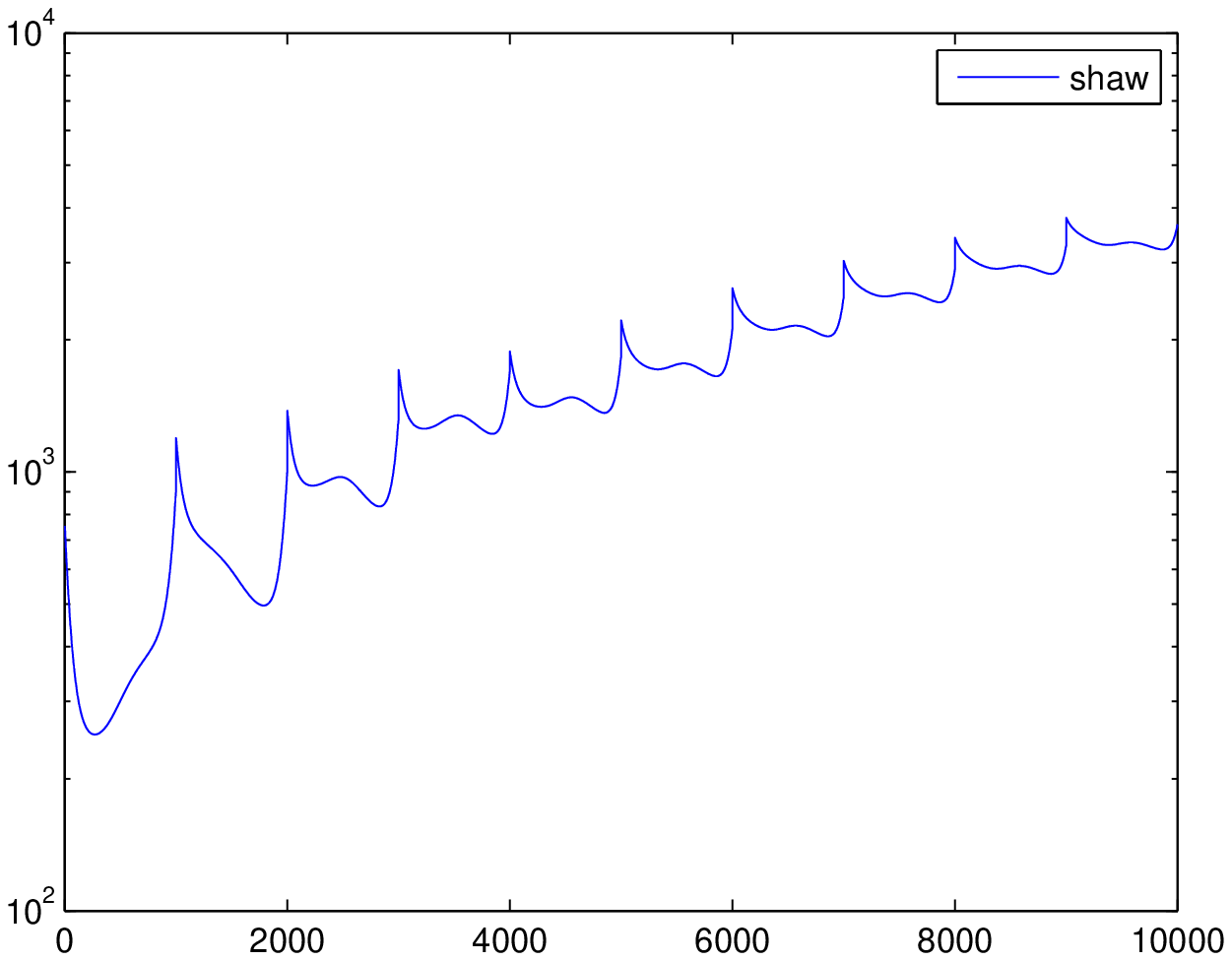}
         \caption{\footnotesize shaw(10000),\protect\\$\delta=0.1$}\label{kacshaw(10000d0.1)}
 \end{minipage}
\end{figure}

\begin{figure}[htb]
  \centering
  \begin{minipage}[h]{.25\linewidth}
      \includegraphics[width=130pt]{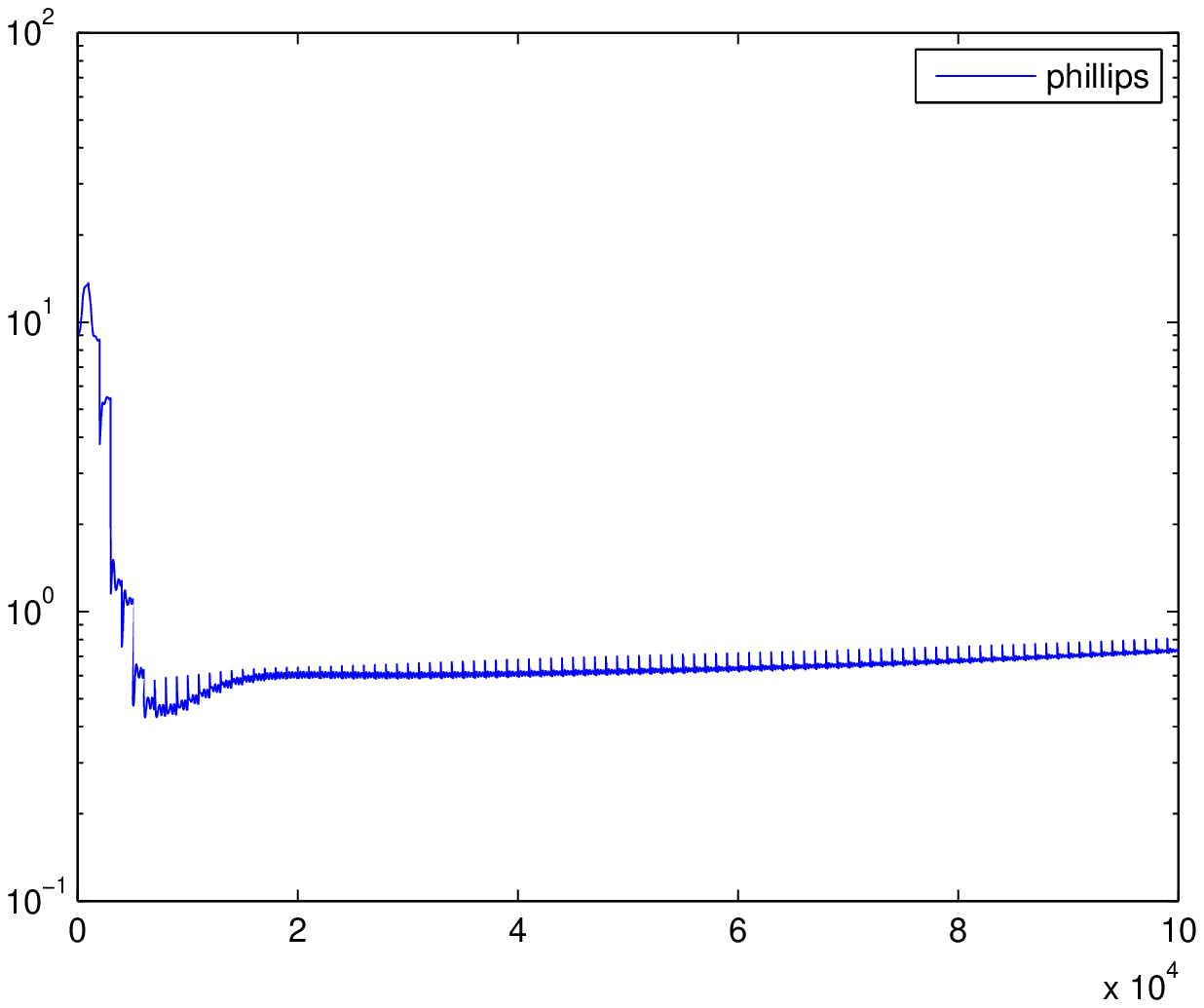}
       \caption{\scriptsize phillips(100000),\protect\\$\delta=0.1$}\label{kacphi(100000d0.1)}
  \end{minipage}
  \begin{minipage}[h]{.25\linewidth}
     \includegraphics[width=130pt]{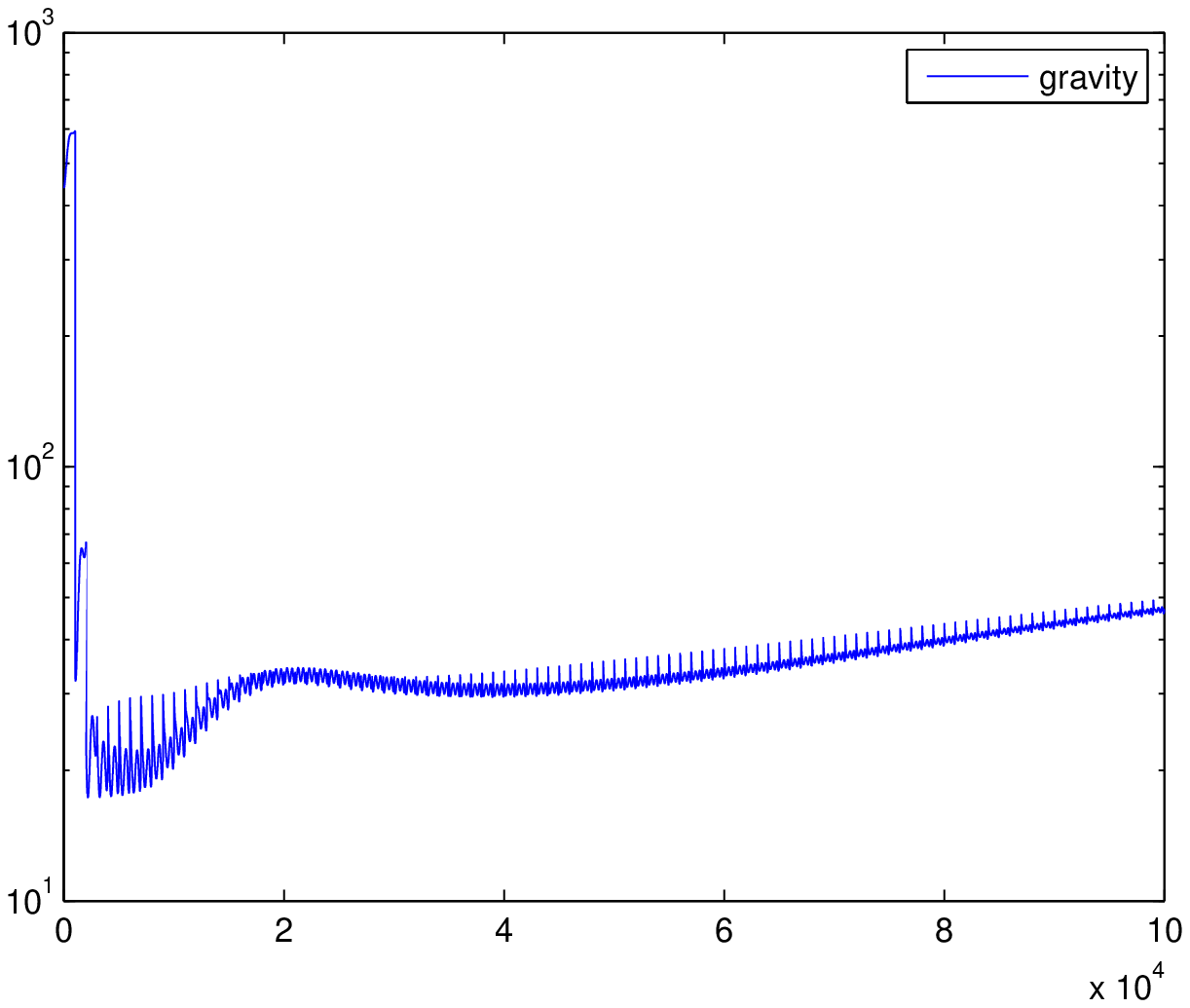}
     \caption{\scriptsize gravity(100000),\protect\\$\delta=0.1$}\label{kacgravity(1000d0.1)}
  \end{minipage}
  \begin{minipage}[h]{.25\linewidth}
     \includegraphics[width=130pt]{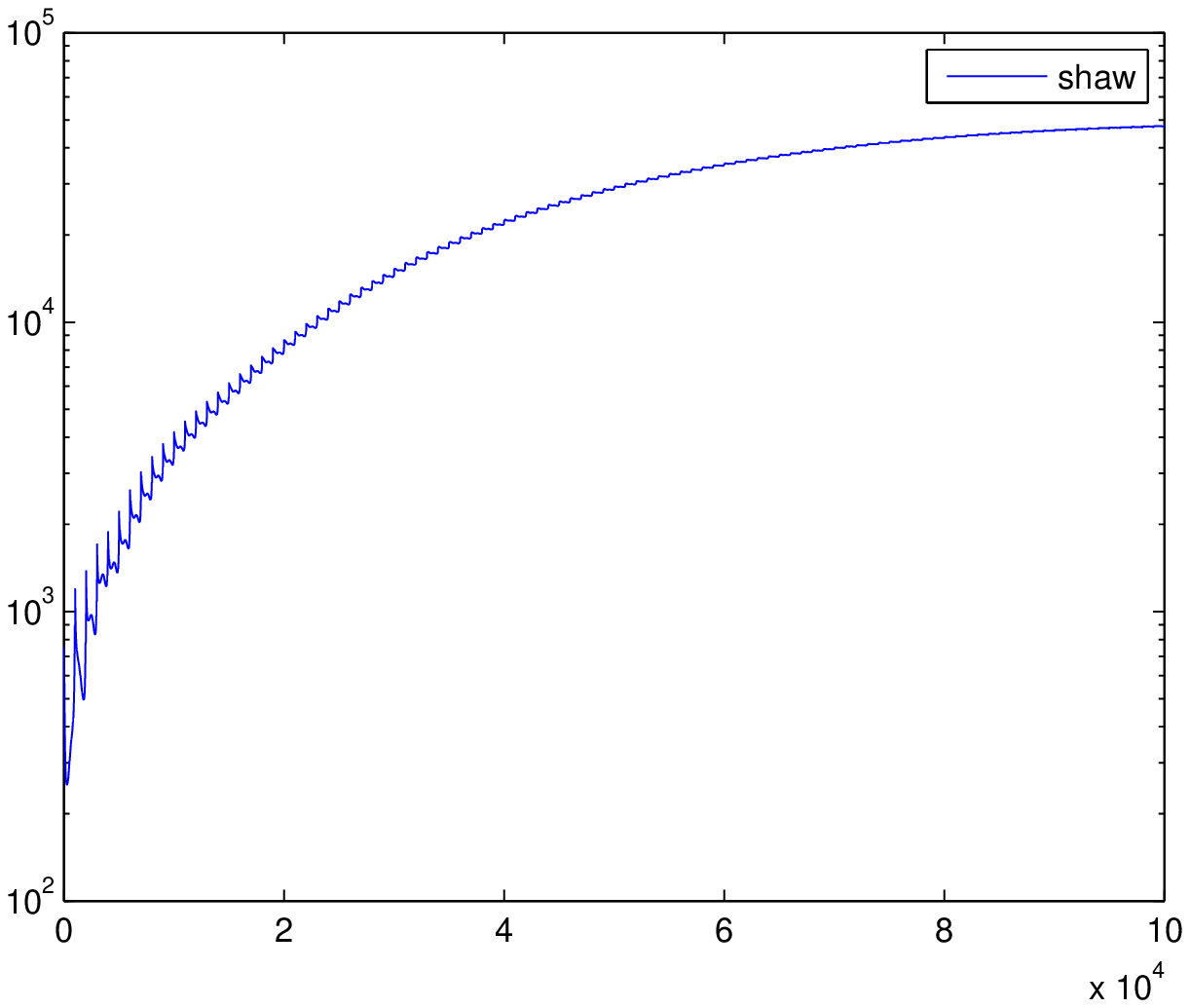}
         \caption{\scriptsize shaw(100000),\protect\\$\delta=0.1$}\label{kacshaw(100000d0.1)}
 \end{minipage}
\end{figure}

Figure \ref{kacphi(10000d0.1)}$\sim$\ref{kacshaw(10000d0.1)} and Figure \ref{kacphi(100000d0.1)}$\sim$\ref{kacshaw(100000d0.1)} show the behavior of the error results for phillips, gravity and shaw solved by Kaczmarz method as $K=10000, \delta=0.1$, respectively. Kaczmarz method takes on 'semi-convergence' in solving phillips, gravity and shaw problems.
\begin{figure}[htb]
  \centering
  \begin{minipage}[h]{.25\linewidth}
      \includegraphics[width=130pt]{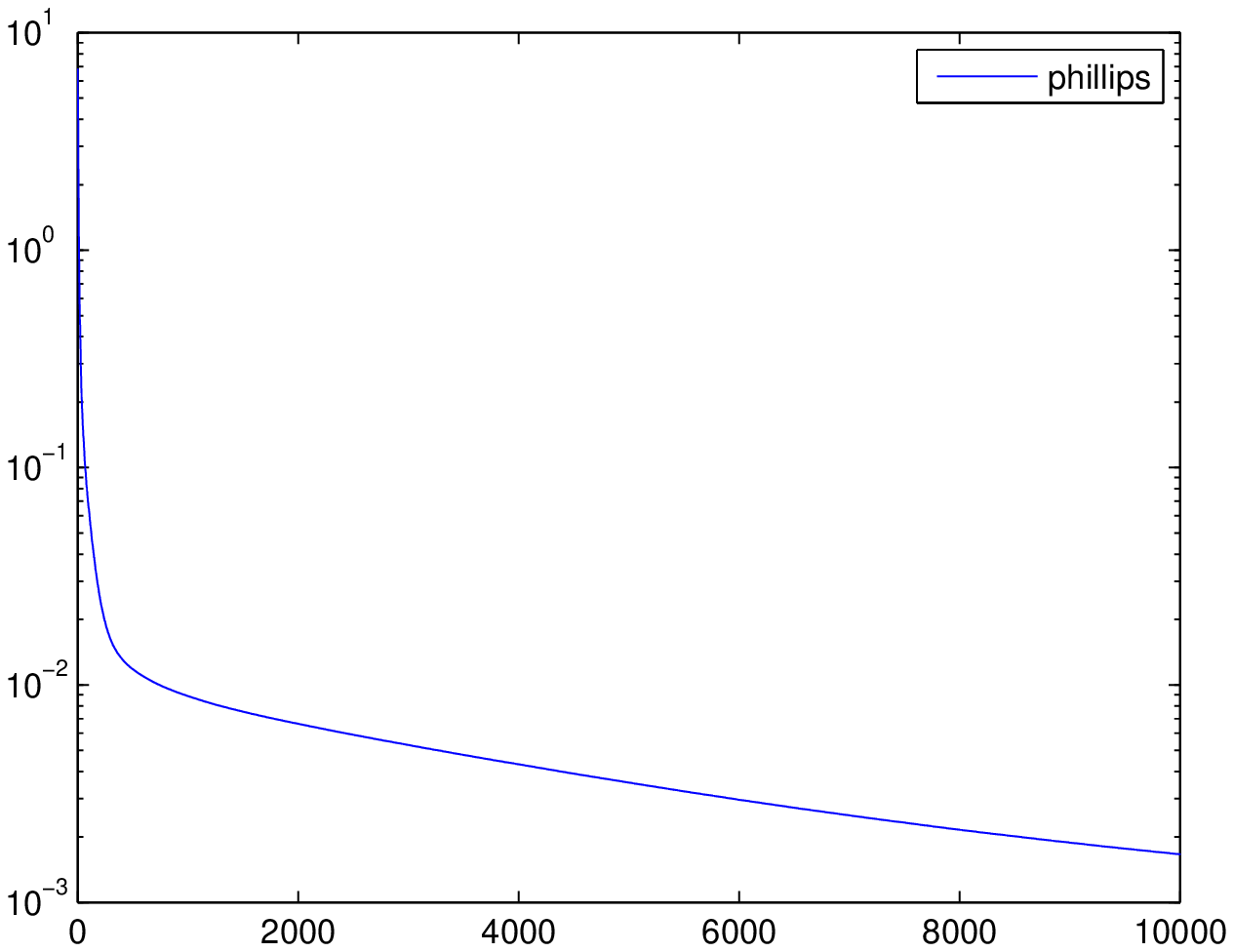}
       \caption{\scriptsize phillips(10000)}\label{kacphi(10000rd0)}
  \end{minipage}
  \begin{minipage}[h]{.25\linewidth}
     \includegraphics[width=130pt]{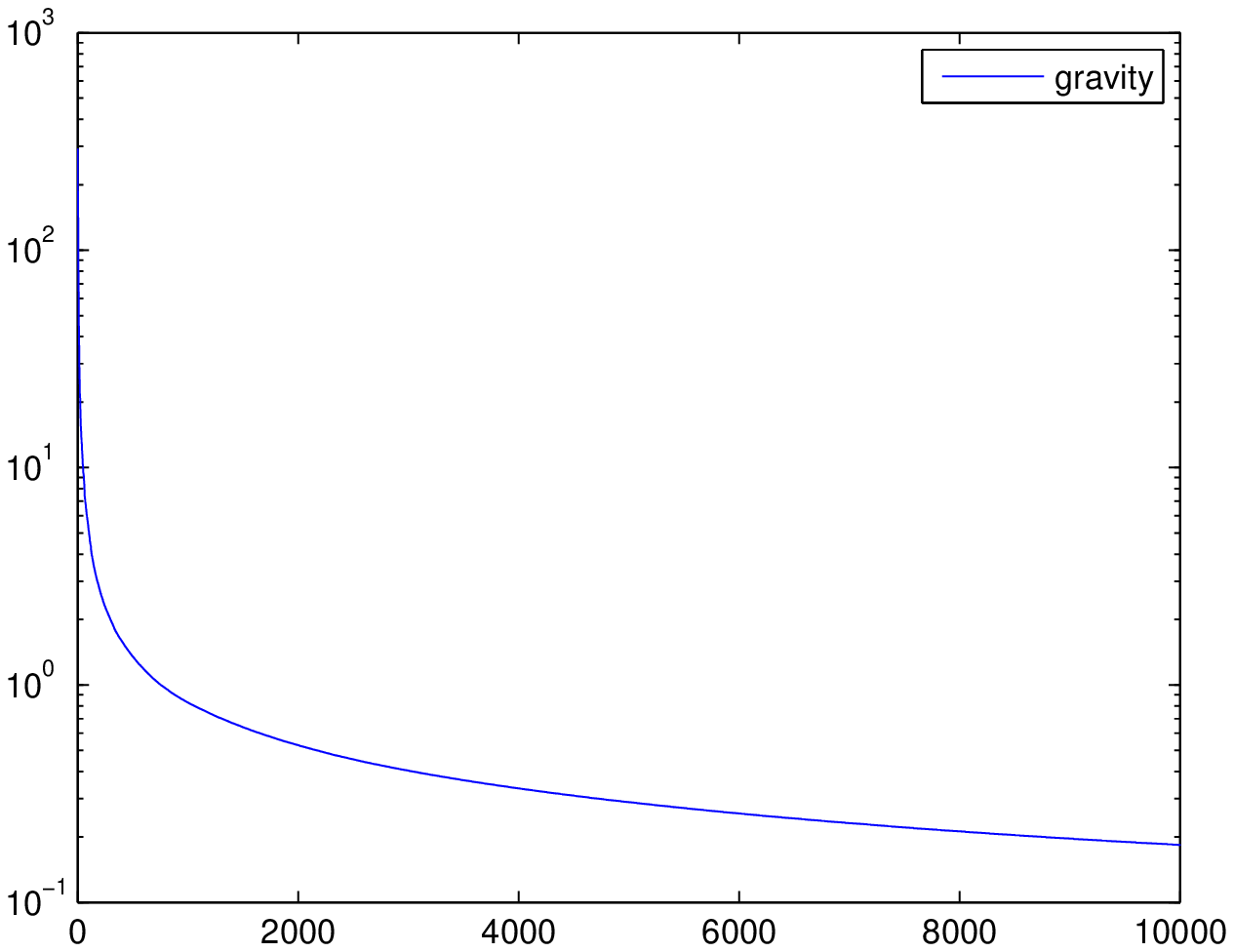}
     \caption{\scriptsize gravity(10000)}\label{kacgravity(10000rd0)}
  \end{minipage}
  \begin{minipage}[h]{.25\linewidth}
     \includegraphics[width=130pt]{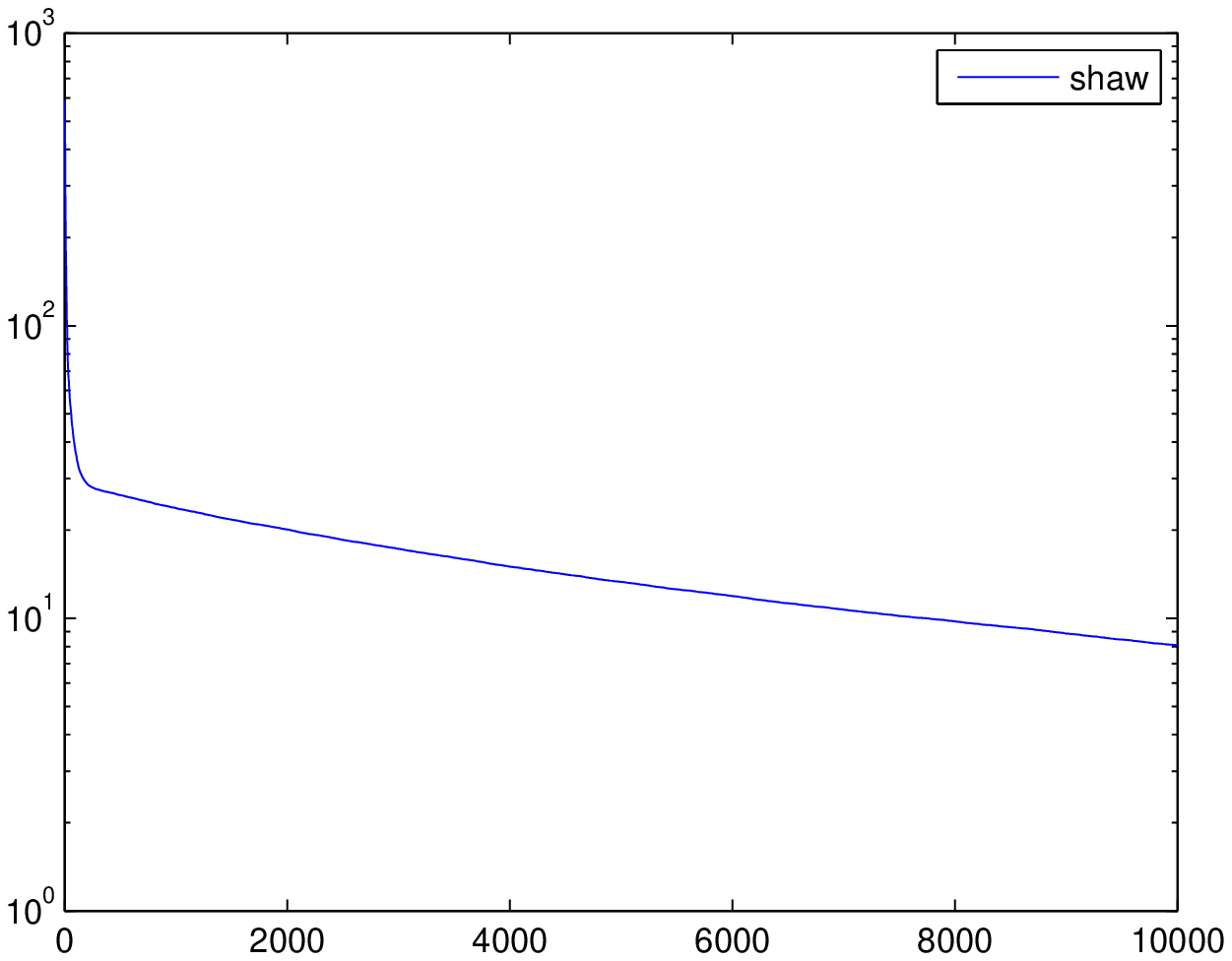}
         \caption{\scriptsize shaw(10000)}\label{kacshaw(10000rd0)}
 \end{minipage}
\end{figure}

\begin{figure}[ht]
  \centering
  \begin{minipage}[h]{.25\linewidth}
      \includegraphics[width=130pt]{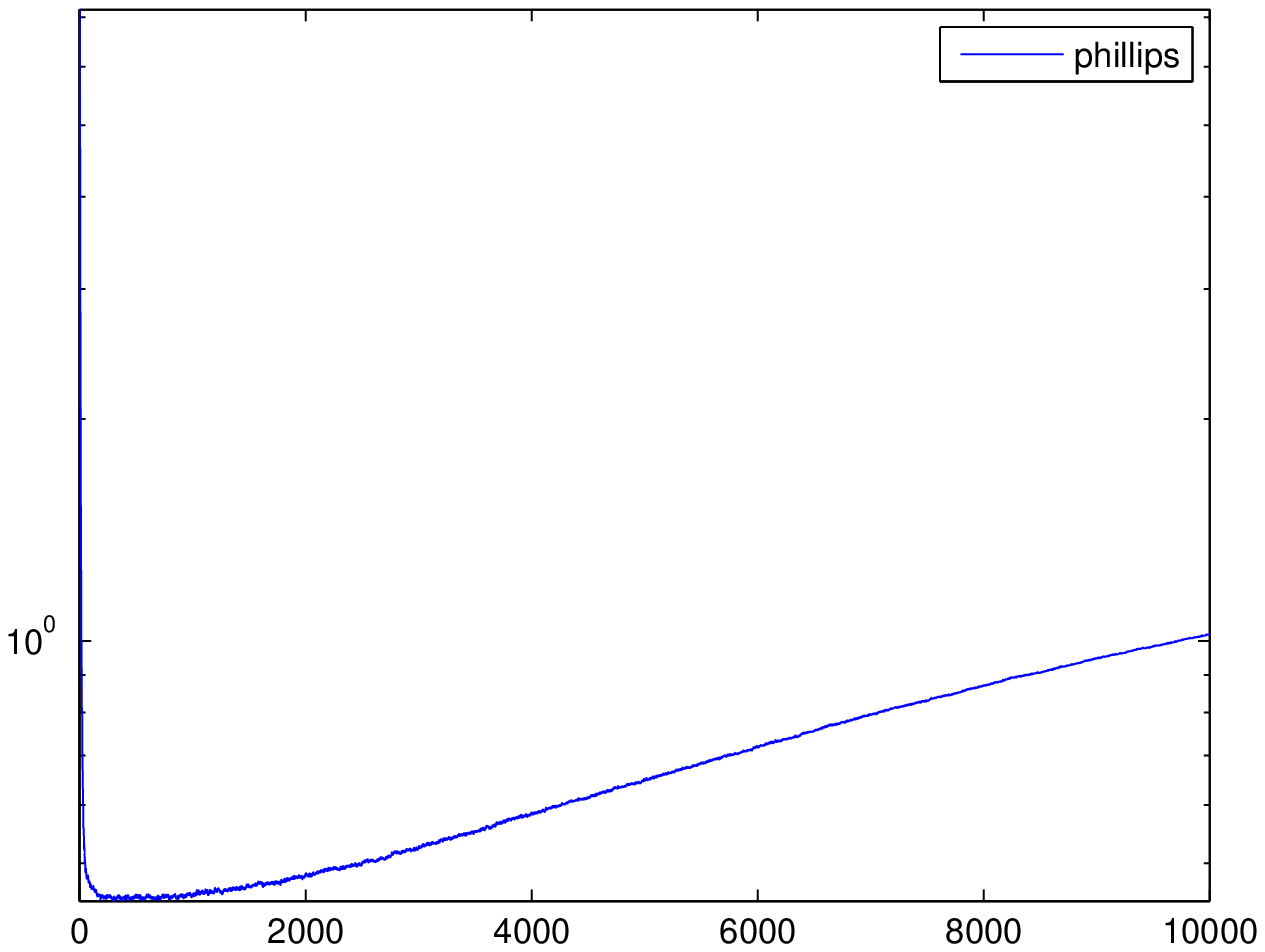}
       \caption{\scriptsize phillips(10000)}\label{kacphi(10000rd0.1)}
  \end{minipage}
  \begin{minipage}[h]{.25\linewidth}
     \includegraphics[width=130pt]{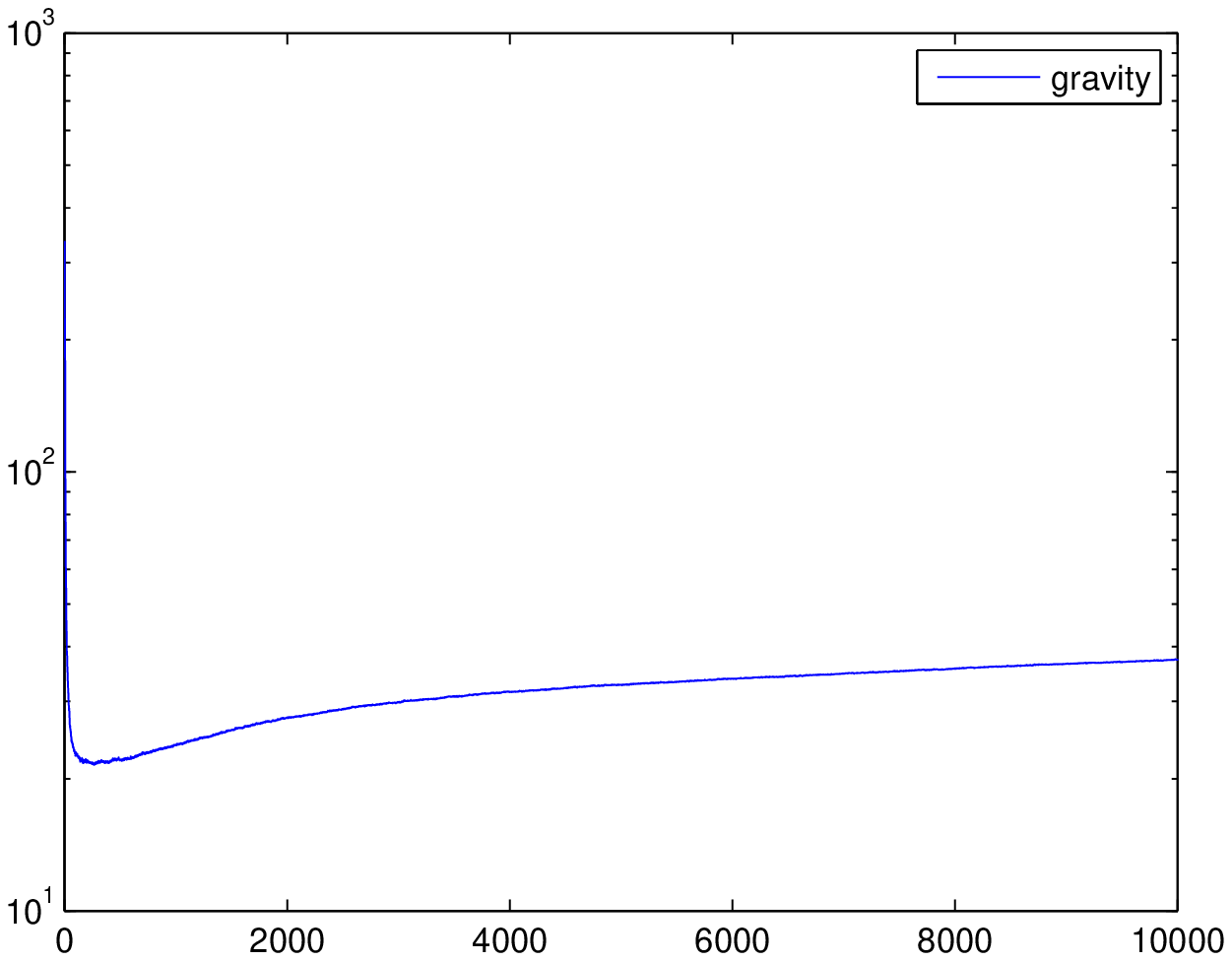}
     \caption{\scriptsize gravity(10000)}\label{kacgravity(10000rd0.1)}
  \end{minipage}
  \begin{minipage}[h]{.25\linewidth}
     \includegraphics[width=130pt]{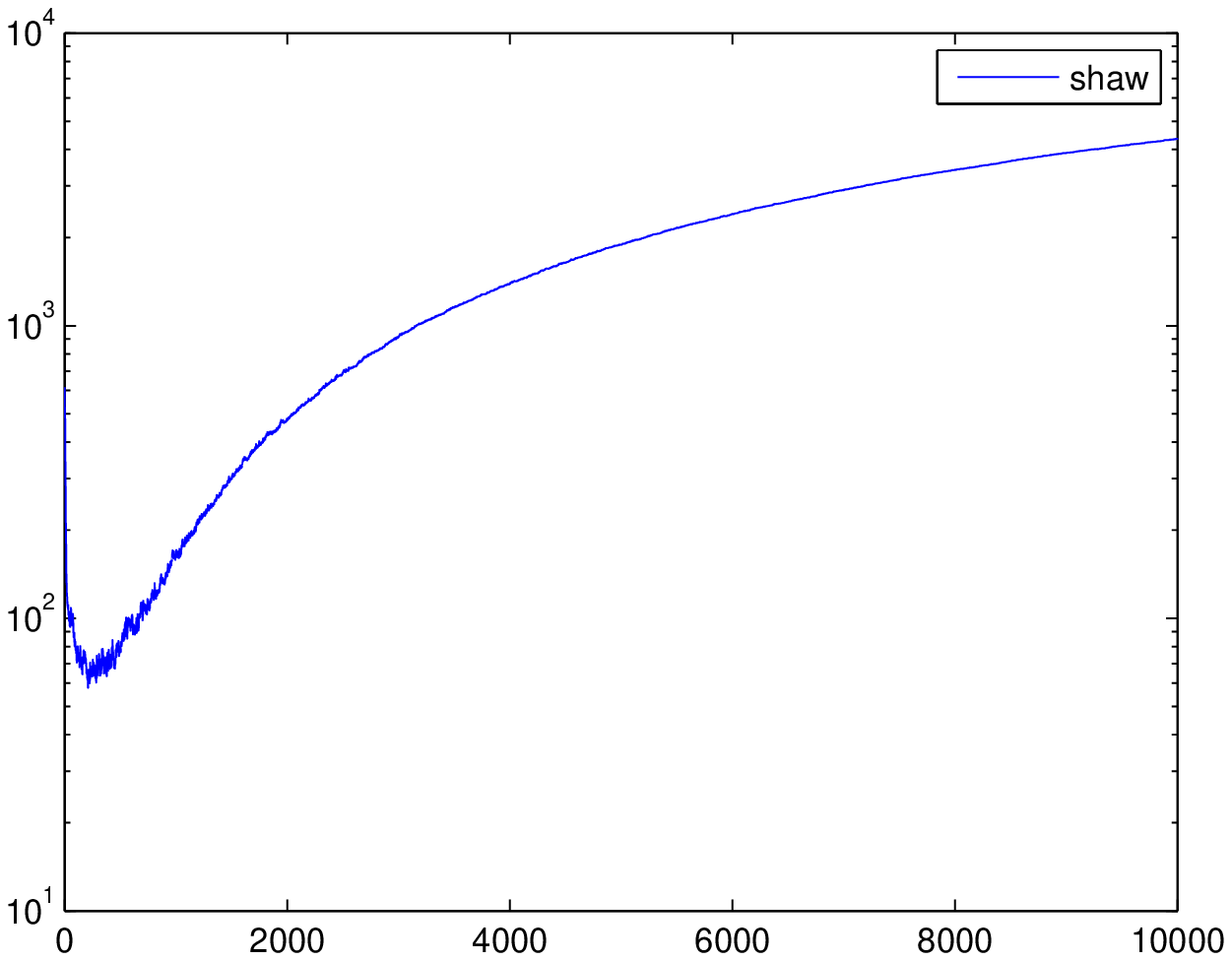}
         \caption{\scriptsize shaw(10000)}\label{kacshaw(10000rd0.1)}
 \end{minipage}
\end{figure}

The behavior of convergence of randomized Kaczmarz mehtod can be seen in Figure \ref{kacphi(10000rd0)}$\sim$ \ref{kacshaw(10000rd0.1)}, where  Figure \ref{kacphi(10000rd0)}$\sim$ \ref {kacshaw(10000rd0)} are the results of randomized Kaczmarz method as $K=10000,\delta=0$, these curve decline monotonously but don't tend to zero, which is to some extent in accordance with Theorem \ref{thm.5}. When there are noise in linear systems, i.e. the right hand of linear system are $b^\delta$, the behavior of convergence of randomized Kaczmarz method are shown in Figure \ref{kacphi(10000rd0.1)}$\sim$ \ref{kacshaw(10000rd0.1)}. Actually, Theorem \ref{thm.5} in this paper couldn't exhibit the tendency of the error on iterative step $k$ as $\delta=0.1$,
nevertheless, Theorem \ref{thm.5} shows the error bound between numerical solution $x_k$ and Moore-Penrose generalized solution $x^\dagger$, and the result of Theorem \ref{thm.randomized.Kaczmarz} is also about the bound rather than the errors about $k$. Not only Theorem \ref{thm.randomized.Kaczmarz} but also Theorem \ref{thm.5} couldn't reflect the 'semi'-convergent behavior.

\section{Conclusion}

In this paper, we consider Kaczmarz like methods for solving linear systems. For consistent systems $Ax=b$,  Theorem \ref{thm.5}, \ref{thm.6}, \ref{theorem.blocker.randomized.kaczamrz} and Corollary \ref{normal_system.cor}, \ref{corollary.7}, \ref{cor.block.randomized.Kaczmarz}  show that the vector sequence $\{x_k\}$ generated from Kaczmarz like methods converge to $P_{N(A)}x_0+x^\dagger$ exponentially. Meanwhile, for inconsistent linear systems, Theorem \ref{thm.5}, \ref{theorem.blocker.randomized.kaczamrz} and Corollary \ref{normal_system.cor}, \ref{cor.block.randomized.Kaczmarz} show Kaczmarz like methods is not convergent. From ill-posed problem theory, the generalized solution $x^\dagger$ is crucial for linear system $Ax=b$.
Therefore, Theorem \ref{thm.5}, \ref{theorem.blocker.randomized.kaczamrz} and Corollary \ref{normal_system.cor},\ref{cor.block.randomized.Kaczmarz} exhibit the convergent tendency to $P_{N(A)}x_0+x^\dagger$ of Kaczmarz method for solving inconsistent linear systems. In addition, most overdetermined linear system $Ax=b^\delta$ is either inconsistent (i.e. $b^\delta=b\notin R(A)$) or ill-posed($b^\delta\neq b$, $b\in R(A)$ and $b^\delta$ is in or not in the range $R(A)$). Hence, Kaczmarz like methods can be regarded as a regularized method, where iterative step $k$ is regularized parameter.

\end{document}